\documentclass[a4paper]{article}
\usepackage{amsmath}
\usepackage{amsfonts}
\usepackage{amssymb}
\usepackage{amsthm}

\def\a{\alpha}
\def\b{\beta}
\def\d{\delta}
\def\e{\epsilon}
\def\g{\gamma}

\def\s{\sigma}

\def\D{\Delta}

\def\X{\bar{X}}
\def\Y{\bar{Y}}
\def\Z{\bar{Z}} 

\def\S{{\tilde{S}}}
\def\ta{{\tilde{\alpha}}}
\def\tb{{\tilde{\beta}}}

\def\ot{\otimes}

\def\R{{\cal R}}

\def\no{\nonumber}

\newcommand{\reff}[1]{(\ref{#1})}
\newcommand{\sect}[1]{\setcounter{equation}{0}\section{#1}}

\newtheorem{Theorem}{Theorem}
\newtheorem{Definition}{Definition}
\newtheorem{Proposition}{Proposition}
\newtheorem{Lemma}{Lemma}

\newtheorem*{Corollary*}{Corollary}

\begin{document}

\begin{titlepage}
\vskip.3in

\begin{center}
{\Large \bf Some Twisted Results}
\vskip.3in
{\large M.D. Gould and T. Lekatsas } 
\vskip.2in
{\em Department of Mathematics, The University of Queensland,\\ Brisbane, Qld 4072, Australia

Email: tel@maths.uq.edu.au
\vskip.2in
\today}
\end{center}

\vskip 2cm

\begin{abstract}
The Drinfeld twist for the opposite quasi-Hopf algebra, $H^{cop}$ is determined and
is shown to be related to the (second) Drinfeld twist on a quasi-Hopf algebra.
The twisted form of the Drinfeld twist is investigated. In the quasi-triangular
case it is shown that the Drinfeld $u$ operator arises from the equivalence of $H^{cop}$
to the quasi-Hopf algebra induced by twisting $H$ with the $R$-matrix. The Altschuler-Coste
$u$ operator arises in a similar way and is shown to be closely related to the Drinfeld
$u$ operator. The quasi-cocycle condition is introduced, and is shown to play a central
role in the uniqueness of twisted structures on quasi-Hopf algebras. A generalisation of
the dynamical quantum Yang-Baxter equation, called the quasi-dynamical quantum Yang-Baxter
equation is introduced.
\end{abstract}

\end{titlepage}
\vskip 3cm


\sect{Introduction\label{intro}}

Quasi-Hopf algebras (QHA) were introduced by Drinfeld~\cite{Dri90} as generalisations of Hopf algebras. QHA
are the underlying algebraic structures of elliptic quantum groups~\cite{Enr97,Fod94,Fel95,Fro97,Jim97,Zha98} and 
hence have an important role in obtaining solutions to the dynamical Yang-Baxter equation. They arise in conformal field
theory~\cite{DPR,DVVV}, algebraic number theory~\cite{DriQ} and in the theory of knots~\cite{ac,Lin00,Mac91}.

The antipode $S$ of a Hopf algebra $H$ is uniquely determined as the inverse of the identity map on $H$ under the convolution product.
For a quasi-Hopf algebra, the triple $(S,\a,\b)$ consisting of the antipode $S$ and canonical elements $\a,\b \in H$ is termed the
{\em quasi-antipode}. The quasi-antipode of a QHA is not unique~\cite{CP,Dri90,Maj}. However, given two QHA which differ only in
their quasi-antipodes, there exists a unique invertible element $v \in H$ relating them. Moreover, to each invertible
element $v \in H$ there corresponds a quasi-antipode, so that the invertible elements $v \in H$ are in bijection with 
the quasi-antipodes. This allows us to work with a fixed choice for the quasi-antipode (more precisely, a fixed equivalence
class for the quasi-antipode). We show that the operator $v \in H$ is universal i.e. invariant under an arbitrary twist
$F \in H \ot H$. In the quasi-triangular case, the equivalence of the quasi-antipode of the opposite QHA
$H^{cop}$ and the quasi-antipode induced by twisting $H$ with the $R$-matrix, gives rise to a specific form of the $v$ operator, which 
we call the Drinfeld - Reshetikhin~\cite{Dri90u,r} $u$ operator. The $u$-operator introduced by Altschuler and Coste~\cite{ac},
arises in a similar way and is shown to be simply related to the Drinfeld - Reshetikhin $u$ operator. In view of the invariance
of the $v$ operators these $u$ operators are also invariant under twisting.

For a Hopf algebra $H$ the antipode $S$ is both an algebra and a co-algebra anti-homomorphism. In the QHA case Drinfeld has shown that
the antipode $S$ is a co-algebra anti-homomorphism only upto conjugation by a twist, $F_\d$ (the Drinfeld twist). Assuming the antipode $S$
is invertible with inverse $S^{-1}$, we show that $S^{-1}$ is a co-algebra anti-homomorphism upto conjugation by an invertible
element $F_0$, which we call the second Drinfeld twist on $H$. The form of the Drinfeld twist for the opposite QHA $H^{cop}$ is determined
and shown to be simply related to this second Drinfeld twist. The behaviour of the Drinfeld twist $F_\d$ under an arbitrary twist
$G \in H \ot H$ is also investigated.

The set of twists on a QHA $H$ form a group. We study a sub-group of the group of twists on a QHA, namely those that
leave the co-product $\D: H \rightarrow H \ot H$ and co-associator $\Phi \in H \ot H \ot H$ unchanged. These twists are
called {\em compatible twists}. Twists that leave the coassociator $\Phi$ unchanged are said to satisfy the quasi-cocycle condition. 
The quasi-cocycle condition is intimately related to the uniqueness of the structure obtained by twisting the quasi-bialgebra
part of a QHA. In the quasi-triangular case we show that $\R^T \R$ and its powers are compatible twists.

Following on from our considerations of the quasi-cocycle condition we introduce the shifted quasi-cocycle
condition on a twist $F(\lambda) \in H \ot H$, where $\lambda \in H$ depends on one (or more) parameters.
We conclude with the quasi-dynamical quantum Yang-Baxter equation (QQYBE), which is the quasi-Hopf
analogue of the usual dynamical QYBE.


\sect{Preliminaries}\label{prelim}
We begin by recalling the definition~\cite{Dri90} of a quasi-bialgebra. 

\begin{Definition}\label{QBA defn}
A quasi-bialgebra $(H,\D,\e, \Phi)$ is a unital associative algebra $H$ over a field $F$,
equipped with algebra homomorphisms $\e: H \rightarrow F$ (co-unit),
$\D: H \rightarrow H \ot H$ (co-product) and an invertible element $\Phi \in H \ot H \ot H$
(co-associator) satisfying
\begin{align}
& (1\ot\D)\D(a)=\Phi^{-1}(\D\ot 1)\D(a)\Phi,\quad \forall a\in H, \label{q-co}\\
& (\D\ot 1\ot 1)\Phi\cdot(1\ot 1\ot\D)\Phi=(\Phi\ot 1)\cdot(1\ot\D\ot 1)
\Phi\cdot(1\ot\Phi), \label{pentagon}\\ 
& (\e\ot 1)\D=1=(1\ot\e)\D, \label{counit}\\ 
& (1\ot\e\ot 1)\Phi=1. \label{epsphi}
\end{align}
\end{Definition}
It follows from equations \reff{pentagon}, \reff{counit} and \reff{epsphi} that the
co-associator $\Phi$ has the additional properties
\begin{align*}
(\e\ot 1 \ot 1)\Phi = 1 = (1 \ot 1 \ot \e) \Phi.
\end{align*}
We now fix the notation to be used throughout the paper.
For the co-associator we follow the notation of \cite{cas,oqhsa} and write
\begin{eqnarray*}
\Phi=\sum_\nu X_\nu\ot Y_\nu\ot Z_\nu~,\qquad 
\Phi^{-1}=\sum_\nu \X_\nu\ot\Y_\nu\ot\Z_\nu.
\end{eqnarray*}
We adopt Sweedler's~\cite{Swe} notation for the co-product 
\begin{eqnarray*}
\D(a)=\sum_{(a)} a_{(1)}\ot a_{(2)}, \qquad\forall a\in H
\end{eqnarray*}
throughout. 
Since the co-product is quasi-coassociative we use the
following extension of Sweedler's notation
\begin{eqnarray}
(1\ot\D)\D(a)&=& a_{(1)}\ot\D(a_{(2)})= a_{(1)}\ot a_{(2)}
^{(1)}\ot a_{(2)}^{(2)} \no \\
(\D\ot 1)\D(a)&=&\D(a_{(1)})\ot a_{(2)}= a_{(1)}^{(1)}
\ot a_{(1)}^{(2)}\ot a_{(2)}. \label{E1}
\end{eqnarray}
In general, the summation sign is omitted from expressions, with the
convention that repeated indices are to be summed over. 

\begin{Definition}\label{QHA defn}
A quasi-Hopf algebra $(H,\D,\e,\Phi,S,\a,\b)$ is a quasi-bialgebra $(H,\D,\e,\Phi)$
equipped with an algebra anti-homomorphism $S$ (antipode) and canonical
elements $\a, \b \in H$ such that
\begin{align}
&  S(X_\nu)\a Y_\nu \b S(Z_\nu)=1=  \X_\nu \b S(\Y_\nu) \a \Z_\nu, \label{Sphi}\\
&  S(a_{(1)})\a a_{(2)} =\e(a) \a, \quad  a_{(1)}\b S(a_{(2)}) =\e(a) \b, 
\quad \forall a\in H \label{Sab}. 
\end{align}
\end{Definition}
Throughout we assume bijectivity of the antipode $S$ so that $S^{-1}$ exists.
The antipode equations~\reff{Sphi},~\reff{Sab} imply $\e(\a)\cdot\e(\b)=1$ and 
$\e(S(a))=\e(S^{-1}(a))=\e(a),\forall a \in H$.
A triple $(S,\a,\b)$ satisfying equations \reff{Sphi},~\reff{Sab} is called a {\em quasi-antipode}.

We shall need the following relations:
\begin{eqnarray}
&& X_\nu a\ot Y_\nu\b S(Z_\nu)= a_{(1)}^{(1)} X_\nu\ot a_{(1)}
^{(2)} Y_\nu\b S(Z_\nu)S(a_{(2)}),\quad \forall a\in H 
\label{aleft}\\
&&\Phi\ot 1\stackrel{\reff{pentagon}}{=}(\D\ot 1\ot1)\Phi\cdot(1\ot1\ot\D)\Phi\cdot(1\ot\Phi^{-1})
\cdot(1\ot\D\ot1)\Phi^{-1} \no \\
&&= X_\nu^{(1)}X_\mu \X_\rho \ot X_\nu^{(2)}Y_\mu\X_\sigma\Y_\rho^{(1)}
\ot Y_\nu Z_\mu^{(1)}\Y_\sigma\Y_\rho^{(2)}\ot Z_\nu Z_\mu^{(2)}\Z_\sigma \Z_\rho \label{phiby1} \\
&&1\ot\Phi=(1\ot\D\ot 1)\Phi^{-1}\cdot(\Phi^{-1}\ot 1)\cdot(\D\ot 1\ot 1)\Phi
\cdot(1\ot 1\ot \D)\Phi \no \\
&&= \X_\nu\X_\mu X_\rho^{(1)}X_\sigma \ot
\Y_\nu^{(1)}\Y_\mu X_\rho^{(2)}Y_\sigma \ot \Y_\nu^{(2)}\Z_\mu Y_\rho Z_\sigma^{(1)}
\ot \Z_\nu Z_\rho Z_\sigma^{(2)} \label{1byphi}
\end{eqnarray}
where we have adopted the notation of equation (\ref{E1}) in \reff{aleft} and the obvious
notation in \reff{phiby1}, \reff{1byphi} so that, for example
\begin{equation*}
\D(X_\nu)= X_\nu^{(1)} \ot X_\nu^{(2)},\quad \mathrm{etc.}
\end{equation*}
Equation~\reff{aleft} follows from applying $(1 \ot m)(1 \ot 1 \ot \b S)$ to
equation~\reff{q-co}, then using~\reff{Sab}.

\section{Uniqueness of the quasi-antipode.}\label{v operator}

For Hopf algebras the antipode $S$ is uniquely determined as the inverse of the identity map
on $H$ under the convolution product. The quasi-antipode $(S,\a,\b)$ for a QHA, is not
unique. Nevertheless it is almost unique as the following result
due to Drinfeld~\cite{Dri90} (whose proof is similar to the one given below) shows:
\begin{Theorem}\label{T1}Suppose $H$ is also a QHA, but with quasi-antipode 
$(\S,\tilde \a ,\tilde \b)$ satisfying \reff{Sphi},\reff{Sab}. Then there
exists a unique invertible $v\in H$ such that
\begin{equation}
v\a=\tilde \a~~,\tilde \b v=\b~~,\S (a)=v S(a)v^{-1}~~,\forall a\in H.\label{E4}
\end{equation}
Explicitly
\begin{eqnarray}
(i)& v= \S(X_\nu)\tilde\a Y_\nu\b S{(Z_\nu)}=\S
{(S^{-1}(\X_\nu))}\S {(S^{-1}(\b))}\S{(\Y_\nu)}\tilde\a\Z_\nu\no\\
(ii)& v^{-1}= S{(X_\nu)}\a Y_\nu\tilde\b\S {(Z_\nu)}=
\X_\nu\tilde\b\S{(\Y_\nu)}\S{(S^{-1}(\a))} \S
{(S^{-1}{(\Z_\nu)})}\no \\
\label{E5}
\end{eqnarray}
\end{Theorem}
\begin{proof}
We proceed stepwise. 

Applying $m\cdot{(\S \ot 1)}{(1 \ot\tilde\a)}$ to equation
\reff{aleft} gives
\begin{equation*}
\S{(X_\nu a)}\ta Y_\nu\b S{(Z_\nu)}=\S{(a_{(1)}^{(1)} X_\nu)}
\ta a_{(1)}^{(2)} Y_\nu\b S{(Z_\nu)}S{(a_{(2)})}
\end{equation*}
so that,
\begin{align}
\S (a)v&=\S{(X_\nu)}\S{(a_{(1)}^{(1)})}\ta a_{(1)}^{(2)} Y_\nu\b
S{(Z_\nu)}S{(a_{(2)})}
\stackrel{\reff{Sab}}{=}vS(a),\quad \forall a\in H \label{S1}
\end{align}
where $m:H \ot H \rightarrow H$ is the multiplication map $m(a \ot b) =ab, \forall a,b \in H$.

Next observe, from equation \reff{phiby1} that, in view of \reff{Sab}, 
\begin{eqnarray}
v\ot 1&=&\S{(X_\nu^{(1)} X_\mu\X_\rho)}\ta X_\nu^{(2)} Y_\mu\X_\s\Y_\rho
^{(1)}\b S{(Y_\nu Z_\mu^{(1)}\Y_\s\Y_\rho^{(2)})}\ot Z_\nu Z_\mu^{(2)}
\Z_\s\Z_\rho\no\\
&=&\S{(X_\mu)}\ta Y_\mu\X_\s\b S{(Z_\mu^{(1)}\Y_\s)}\ot Z_\mu^{(2)}\Z_\s. \no
\end{eqnarray}
Applying $m\cdot(1 \ot \a)$ from the left gives
\begin{align}
v\a=&\S{(X_\mu)}\ta Y_\mu\X_\s\b S{(Z_\mu^{(1)}\Y_\s)}\a Z_\mu^{(2)} \Z_\s \no \\
=&\ta\X_\s\b S{(\Y_\s)}\a\Z_\s\stackrel{\reff{Sphi}}{=}\ta. \label{S2}
\end{align}
From this it follows that
\begin{eqnarray}
&&\S{(S^{-1}{(\X_\nu)})}\cdot\S{(S^{-1}(\b))}\cdot\S{(\Y_\nu)}\ta\Z_\nu \no\\
&\stackrel{(\ref{S2})}{=}&\S{(S^{-1}{(\X_\nu)})}\cdot\S{(S^{-1}(\b))}\S{(\Y_\nu)}\cdot v\a\Z_\nu\no\\
&\stackrel{(\ref{S1})}{=}&v\cdot S{(S^{-1}{(\X_\nu)})}\cdot S{(S^{-1}(\b))}\cdot S{(\Y_\nu)}\a\Z_\nu\no\\
&=&v\cdot \X_\nu\b S{(\Y_\nu)}\a\Z_\nu \stackrel{\ref{Sphi}}{=}v \no
\end{eqnarray}
which proves (\ref{E5})(i). To see $v$ is invertible observe that
\begin{eqnarray}
v\cdot  S{(X_\nu)}\a Y_\nu\tb\S{(Z_\nu)}&\stackrel{(\ref{S1})}{=}&\S{(X_\nu)}v\a Y_\nu\tb\S{(Z_\nu)} \no \\
&\stackrel{(\ref{S2})}{=}&\S{(X_\nu)}\ta Y_\nu\tb\S{(Z_\nu)} \no \\
&\stackrel{\reff{Sphi}}{=}&1 \no
\end{eqnarray}
so
\begin{equation*}
v^{-1}= S{(X_\nu)}\a Y_\nu\tb \S{(Z_\nu)}
\end{equation*}
as stated.

Now using equation \reff{1byphi} we have
\begin{align}
1\ot v^{-1}&= \X_\nu\X_\mu X_\rho^{(1)} X_\s\ot S{(\Y_\nu^{(1)}\Y_\mu X_\rho^{(2)}
Y_\s)}\a \Y_\nu^{(2)}\Z_\mu Y_\rho\Z_\s^{(1)}\tb\S{(\Z_\nu Z_\rho Z_\s^{(2)})}\no\\
&\stackrel{\reff{Sab}}{=}\X_\mu X_\rho^{(1)}\ot S{(\Y_\mu X_\rho^{(2)})}\a\Z_\mu Y_\rho\tb\S{(Z_\rho)}.\no
\end{align}
Applying $m \cdot (1 \ot \b)$ gives
\begin{align}
\b v^{-1}=&\X_\mu X_\rho^{(1)}\b S{(\Y_\mu X_\rho^{(2)})}\a\Z_\mu Y_\rho
\tb\S{(Z_\rho)}\no\\
=&\X_\mu\b S{(\Y_\mu)}\a\Z_\mu\cdot\tb\stackrel{\reff{Sphi}}{=}\tb \label{S3}
\end{align}
which completes the proof of (\ref{E4}). As to (\ref{E5})(ii) observe that
\begin{eqnarray}
&&\tilde{X}_\nu\tb\S{(\Y_\nu)}\S{(S^{-1}(\a))}\S{(S^{-1}{(\Z_\nu)})}\no\\
&\stackrel{(\ref{S3})}{=}&\X_\nu \b v^{-1}\S{(\Y_\nu)}\S{(S^{-1}(\a))}\S{(S^{-1}{(\Z_\nu)})}\no\\
&\stackrel{(\ref{S1})}{=}&\X_\nu\b S{(\Y_\nu)}S{(S^{-1}(\a))}S{(S^{-1}{(\Z_\nu)})}v^{-1}\no\\
&=&\X_\nu\b S{(\Y_\nu)}\a\Z_\nu\cdot v^{-1}\stackrel{\reff{Sphi}}{=}v^{-1}\no
\end{eqnarray}
as required. It finally remains to prove uniqueness. Hence suppose $u\in H$
satisfies
\begin{equation*}
u S (a)=\S(a)u,\quad \forall a\in H,~~u\a=\ta,~~\tb u=\b.
\end{equation*}
Then
\begin{eqnarray}
uv^{-1}&=&u\cdot S{(X_\nu)}\a Y_\nu\tb\S{(Z_\nu)}\no\\
&=& \S{(X_\nu)}u\a Y_\nu\tb\S{(Z_\nu)}\no\\
&=& \S{(X_\nu)}\ta Y_\nu\tb\S{(Z_\nu)}\stackrel{\reff{Sphi}}{=}1\no
\end{eqnarray}
which implies $u=v$ as required.
\end{proof}

In the special case $\S=S$ we obtain the following useful result.
\begin{Corollary*}
Suppose $H$ is also a QHA with quasi-antipode $(S,\ta,\tb)$. Then there is a unique
invertible central element $v\in H$, given explicitly by equation (\ref{E5})(i)
(with $\S=S$), such that
\begin{equation*}
v\a=\ta~~~~,~~~~\tb v=\b.
\end{equation*}
\end{Corollary*}
It thus follows that the triple $(S,\a,\b)$ satisfying \reff{Sphi}, \reff{Sab} for a QHA
is not unique. Indeed following Theorem~\ref{T1}, for arbitrary invertible
$v\in H$, the triple $(\S,\ta,\tb)$ defined by
\begin{equation*}
\S(a)=v S(a)v^{-1},\quad \forall a\in H~~~;\ta=v\a~~~,\tb=\b v^{-1}
\end{equation*}
is easily seen to satisfy \reff{Sphi}, \reff{Sab} and thus gives rise to a 
quasi-antipode $(\S,\ta,\tb)$. Theorem~\ref{T1}
then shows that all such quasi-antipodes $(\S,\ta,\tb)$ are
obtainable this way: thus there is a 1--1 correspondence between the latter
and invertible $v\in H$. We say that these structures are {\it equivalent} since they 
clearly give rise to equivalent QHA structures. Throughout we work with a fixed
choice for the quasi-antipode $(S,\a,\b)$.

We conclude this section with the following useful result, proved in
\cite{oqhsa}, concerning the opposite QHA structure on $H$:
\begin{Proposition}\label{P1}
$H$ is also a QHA, with co-unit $\e$, under the opposite co-product and co-associator
$\D^T$, $\Phi^T \equiv \Phi_{321}^{-1}$ respectively, with quasi-antipode
$(S^{-1},\a^T=S^{-1}(\a),~\b^T=S^{-1}(\b))$.
\end{Proposition}
The QHA $H^{cop} \equiv (H, \D^T, \e, \Phi^T,S^{-1},\a^T,\b^T)$ is called the opposite QHA structure.
We remark that above we have adopted the notation of~\cite{cas} and~\cite{oqhsa} so that $\D^T
=T\cdot\D,~~T$ the usual twist map, and
\begin{equation*}
\Phi_{321}^{-1}= \Z_\nu\ot\Y_\nu\ot\X_\nu.
\end{equation*}
This latter notation extends in a natural way and will be employed throughout.

\section{Twisting}\label{twisting}

Let $H$ be a quasi-bialgebra. Then $F \in H \ot H$ is called a twist
if it is invertible and satisfies the co-unit property
\begin{equation*}
{(\e \ot 1)}F={(1\ot\e)}F=1.
\end{equation*}
We recall that $H$ is also a QBA with the same co-unit $\e$ but with co-product
and co-associator given by
\begin{eqnarray}
&&\D_F (a)=F\D(a)F^{-1},\quad \forall a\in H\no\\
&&\Phi_F={(F\ot 1)}\cdot{(\D\ot 1)}F\cdot\Phi\cdot{(1 \ot\D)F^{-1}}\cdot
{(1\ot F^{-1})},\label{E6}
\end{eqnarray}
called the twisted structure induced by $F$. If moreover $H$ is a QHA with
quasi-antipode $(S,\a,\b)$ then $H$ is also a QHA under the
above twisted structure with the {\it same} antipode $S$ but with canonical elements
\begin{equation}
\a_F=m\cdot{(1\ot\a)}{(S\ot 1)}F^{-1}\quad\quad\b_F=m\cdot{(1\ot\b)}{(1\ot S)}F\label{E7}
\end{equation}
respectively. A detailed proof of these well known results is given in~\cite{Zha98}.
We now investigate the behaviour of the operator $v$ of Theorem~\ref{T1} under the 
twisted structure induced by $F$.

\subsection{Universality of $v$}\label{Uni of v} 

Recall that the operator $v$ is given by
\begin{equation*}
v=\S(X_\nu)\tilde{\a} Y_\nu \b S(Z_\nu)
\end{equation*}
Let $F \in H\ot H$ be an arbitrary twist. We use the following notation for the twist $F$
and its inverse $F^{-1}$,
\begin{equation*}
F= f_i \ot f^i ,~~~F^{-1}= \bar{f}_i \ot \bar{f}^i.
\end{equation*}
The twisted form of the co-associator is given by~(\ref{E6})
\begin{align}
\Phi_F&= X_\nu^F \ot Y_\nu^F \ot Z_\nu^F =
 f_i f_j^{(1)}X_\nu \bar{f}_k \ot f^i f_j^{(2)}Y_\nu\bar{f}_{(1)}^k\bar{f}_l
\ot f^j Z_\nu \bar{f}^k_{(2)}\bar{f}^l. \label{SS1}
\end{align}
For the twisted forms of the canonical elements we have from~(\ref{E7})
\begin{align}
\tilde{\a}_F&=m \cdot (1 \ot \tilde{\a})(\S \ot 1)F^{-1}=\S(\bar{f}_p)\tilde{\a}\bar{f}^p \no \\
\b_F&=m \cdot (1 \ot \b)(1 \ot S)F= f_q \b S(f^q).  \label{SS2}
\end{align}
We note that
\begin{align}
\S(f_j)\tilde{\a}_F f^j&\stackrel{(\ref{SS2})}{=}\S(\bar{f}_pf_j)\tilde{\a} \bar{f}^p f^j
=m \cdot (1 \ot \a)(\S \ot 1)(F^{-1}F)=\tilde{\a} \label{SS3}
\end{align}
and similarly,
\begin{align}
\bar{f}_j \b_F S(\bar{f}^j)=\b. \label{SS4}
\end{align}
The twisted form of $v$ is given by
\begin{align}
v_F&~~=~~\S(X^F_\nu)\tilde{\a}_F Y^F_\nu \b_F S(Z^F_\nu) \no \\
&\stackrel{(\ref{SS1})}{=} \S(f_if_j^{(1)}X_\nu\bar{f}_k)\tilde{\a}_F f^if_j^{(2)} Y_\nu \bar{f}^k_{(1)} \bar{f}_l
\b_F S(f^j Z_\nu \bar{f}^k_{(2)} \bar{f}^l) \no \\
&~~=~~ \S(f_j^{(1)}X_\nu\bar{f}_k)\S(f_i) \tilde{\a}_F f^if_j^{(2)} Y_\nu \bar{f}^k_{(1)} \bar{f}_l
\b_F S(\bar{f}^l) S(f^j Z_\nu \bar{f}^k_{(2)} ) \no \\
&\stackrel{(\ref{SS3})}{=}\S(f_j^{(1)}X_\nu\bar{f}_k) \tilde{\a} f_j^{(2)} Y_\nu \bar{f}^k_{(1)} \bar{f}_l
\b_F S(\bar{f}^l) S(f^j Z_\nu \bar{f}^k_{(2)} ) \no \\
&\stackrel{(\ref{SS4})}{=}\S(f_j^{(1)}X_\nu\bar{f}_k) \tilde{\a} f_j^{(2)} Y_\nu \bar{f}^k_{(1)} 
\b S(f^j Z_\nu \bar{f}^k_{(2)} ) \no \\
&~~=~~ \S(X_\nu\bar{f}_k) \S(f_j^{(1)}) \tilde{\a} f_j^{(2)} Y_\nu \bar{f}^k_{(1)} 
\b S(\bar{f}^k_{(2)}) S(f^j Z_\nu ) \no \\
&~~=~~ \S(X_\nu\bar{f}_k)  \tilde{\a}  Y_\nu \bar{f}^k_{(1)} 
\b S(\bar{f}^k_{(2)}) S(Z_\nu ) \no \\
&~~=~~ \S(X_\nu)  \tilde{\a}  Y_\nu \b S(Z_\nu ) = v, \no 
\end{align}
where, in the last two lines we have used the antipode properties of $\a,\b$~\reff{Sab} and
the co-unit property of twists.
We have thus proved
\begin{Theorem}\label{uni_of_v}
The operator $v$ is universal (i.e. invariant under twisting).
\end{Theorem}

\section{The Drinfeld twists}\label{Drinfeld twists}
We turn our attention to the Drinfeld twist for the opposite structure of 
proposition~\ref{P1}. It is tempting to assume that $F_\d^T$ qualifies as a Drinfeld twist
for the opposite structure. However this is not true since the antipode for the latter is
$S^{-1}$ rather than $S$. We shall show that the Drinfeld twist for the opposite 
structure is in fact related to the second Drinfeld twist which we define below. We
begin with a review of the Drinfeld twist.

\subsection{The Drinfeld twist}
Observe that $\D'$ defined by
\begin{equation}
\D'(a)={(S\ot S)}\D^T{(S^{-1}(a))},\quad \forall a\in H \label{E8}
\end{equation}
also determines a co-product on $H$. Associated with this co-product we have
a new QHA structure on $H$, which was proved in~\cite{oqhsa} and 
which we restate here:
\begin{Proposition}\label{P2}
$H$ is also a QHA with the same co-unit $\e$ and antipode $S$ but with co-product $\D'$, 
co-associator $\Phi'={(S\ot S\ot S)}\Phi_{321}$, and canonical elements 
$\a'= S(\b),~\b'=S(\a)$ respectively. 
\end{Proposition}

Drinfeld has proved the remarkable result that this QHA structure is
obtained by twisting with the Drinfeld twist, herein denoted $F_\d$, given explicitly by
\begin{eqnarray}
&(i)&F_\d={(S\ot S)}\D^T{(X_\nu)}\cdot\gamma\cdot\D
{(Y_\nu\b S{(Z_\nu)})}\no\\
&\mbox{}&~~~~=\D'{(\X_\nu\b S{(\Y_\nu)})}\cdot\gamma\cdot\D{(\Z_\nu)} \no
\end{eqnarray}
where
\begin{eqnarray}
&(ii)&\gamma= S{(B_i)}\a C_i\ot S{(A_i)}\a D_i \no
\end{eqnarray}
with
\begin{eqnarray}
&(iii)& A_i\ot B_i\ot C_i\ot D_i=\left\{
\begin{array}{l}
(\Phi^{-1}\ot 1) \cdot(\D \ot 1 \ot 1) \Phi \\
\textrm{or} \\
(1 \ot \Phi)\cdot (1 \ot 1 \ot \D)\Phi^{-1}~.
\end{array} 
\right. \no \\ \label{E9}
\end{eqnarray}
The inverse of $F_\d$ is given explicitly by
\begin{eqnarray}
&(i)& F_\d^{-1}=\D{(\X_\nu)}\cdot\bar\gamma\cdot\D'
{(S{(\Y_\nu)}\a\Z_\nu)}\no\\
&\mbox{}&~~~~~~~=\D{(S{(X_\nu)}\a Y_\nu)}\cdot\bar\gamma\cdot{(S\ot S)}\D^T
{(Z_\nu)}\no
\end{eqnarray}
where
\begin{eqnarray}
&(ii)&\bar\gamma=\bar A_i\b S{(\bar D_i)}\ot\bar B_i\b S{(\bar C_i)} \no
\end{eqnarray}
with
\begin{eqnarray}
&(iii)& \bar A_i\ot \bar B_i\ot \bar C_i\ot \bar D_i=\left\{
\begin{array}{l}
(\D \ot 1 \ot 1) \Phi^{-1} \cdot (\Phi\ot 1) \\
\textrm{or} \\
(1 \ot 1 \ot \D)\Phi\cdot (1 \ot \Phi^{-1})~. 
\end{array} 
\right. \no \\ \label{E10}
\end{eqnarray}
The detailed proof that the QHA structure of proposition~\ref{P2} is
obtained by twisting with $F_\d$, as given in (\ref{E9}), and in particular
\begin{equation}
\D'(a)=F_\d\D(a)F_\d^{-1},\quad \forall a\in H \label{E11}
\end{equation}
is proved in~\cite{oqhsa}. We simply state here some properties of $\gamma, \bar\gamma$
proved in~\cite{oqhsa} and which are crucial to the demonstration of Drinfeld's result:
\begin{Proposition}\label{P3} 
\begin{eqnarray}
&(i)& {(S\ot S)}\D^T{(a_{(1)})}\cdot\gamma\cdot\D{(a_{(2)})}=\e(a)\gamma,\quad \forall a \in H \no\\
&(ii)& \D{(a_{(1)})}\cdot\bar\gamma\cdot{(S\ot S)}\D^T{(a_{(2)})}=
\e (a)\bar\gamma,\quad \forall a \in H \no\\
&(iii)& F_\d\D(\a)=\gamma~~~,~~~\D(\b)F_\d^{-1}~~=~~\bar\gamma~.\label{E12}
\end{eqnarray}
\end{Proposition}
\vskip 3mm
\subsection{The second Drinfeld twist}\label{DT2}

Replacing $S$ with $S^{-1}$ we obtain yet another co-product $\D_0$ on $H$:
\begin{equation}
\D_0(a)={(S^{-1}\ot S^{-1})}\D^T{(S(a))},\quad \forall a\in H. \tag{\ref{E8}$'$}
\end{equation}
We have the following analogue of proposition~\ref{P2}, the proof of which
parallels that of~\cite{oqhsa} proposition 4, but with $S$ and $S^{-1}$ interchanged:
\vskip 3mm
\noindent{\bf Proposition \ref{P2}$~'$~~}
{\it $H$ is also a QHA with the same co-unit $\e$ and antipode $S$ but with co-product
$\D_0$, co-associator $\Phi_0={(S^{-1}\ot S^{-1}\ot S^{-1})}\Phi_{321}$ and 
canonical elements $\a_0=S^{-1}(\b),~\b_0=S^{-1}(\a)$ respectively.
}
\vskip 3mm
By symmetry we would expect this structure to be obtainable twisting. 
Indeed we have
\begin{Theorem}\label{T2}: The QHA structure of proposition~\ref{P2}$~'$ is 
obtained by twisting with
\begin{equation}
F_0\equiv{(S^{-1}\ot S^{-1})}F_\d^T\label{E13}
\end{equation}
herein referred to as the second Drinfeld twist,
where $F_\d$ is the Drinfeld twist and $F_\d^T=T\cdot F_\d.$
\end{Theorem}
\begin{proof}
It is clear that $F_0$ is invertible with inverse
$F_0^{-1}={(S^{-1}\ot S^{-1})}{(F_\d^T)}^{-1}$ and qualifies
as a twist. For the co-product we observe,
\begin{eqnarray}
F_0\D(a)F_0^{-1}&=&{(S^{-1}\ot S^{-1})}F_\d^T\cdot\D(a)\cdot
{(S^{-1}\ot S^{-1})}{(F_\d^T)}^{-1}\no\\
&=&{(S^{-1}\ot S^{-1})}\cdot T\cdot{[F_\d^{-1}\cdot{(S\ot S)}
\D^T(a)\cdot F_\d]}\no\\
&=&{(S^{-1}\ot S^{-1})}\cdot T\cdot{[F_\d^{-1}\D'{(S(a))}F_\d]}\no\\
&\stackrel{(\ref{E11})}{=}&{(S^{-1}\ot S^{-1})}\cdot T\cdot\D{(S(a))}={(S^{-1}\ot S^{-1})}
\D^T{(S(a))}\no \\
&\stackrel{(\ref{E8}')}{=}&\D_0(a), \quad \forall a\in H.\no
\end{eqnarray}         
The co-associator is slightly more complicated, though also simple. We
have from Drinfeld's result
\begin{equation*}
\Phi'\equiv {(S\ot S\ot S)}\Phi_{321}={(F_\d\ot 1)}\cdot{(\D\ot 1)}
F_\d\cdot\Phi\cdot{(1\ot\D)}F_\d^{-1}\cdot{(1\ot F_\d^{-1})}
\end{equation*}
which implies
\begin{equation*}
{(S\ot S\ot S)}\Phi={[{(F_\d\ot 1)}\cdot{(\D\ot 1)}F_\d\cdot\Phi
\cdot{(1\ot\D)}F_\d^{-1}\cdot{(1\ot F_\d^{-1})}]}_{321}
\end{equation*}
\begin{equation*}
={(1\ot F_\d^T)}\cdot{(1\ot\D^T)}F_\d^T\cdot\Phi_{321}\cdot
{(\D^T\ot 1)}{(F_\d^T)}^{-1}\cdot{(F_\d^{T-1}\ot 1)}.
\end{equation*}
Applying ${(S^{-1}\ot S^{-1}\ot S^{-1})}$ gives
\begin{eqnarray}
\Phi&=&{(F_0^{-1}\ot 1)}\cdot{(\D_0\ot 1)}F_0^{-1}\cdot\Phi_0\cdot
{(1\ot\D_0)}F_0\cdot{(1\ot F_0)}\no\\
&=&{(\D\ot 1)}F_0^{-1}\cdot{(F_0^{-1}\ot 1)}\cdot\Phi_0\cdot
{(1\ot F_0)}\cdot{(1\ot \D)}F_0\no
\end{eqnarray}
with $F_0$ as in the Theorem. Thus
\begin{equation*}
\Phi_0={(F_0\ot 1)}\cdot{(\D\ot 1)}F_0\cdot\Phi\cdot{(1\ot\D)}
F_0^{-1}\cdot{(1\ot F_0^{-1})}
\end{equation*}
which shows that indeed $\Phi_0$ is obtained from $\Phi$ by twisting
with $F_0$. The proof for the canonical elements is straightforward.
\end{proof}

\subsection{The Drinfeld twists for the opposite structure}

Recall that under the opposite structure of proposition~\ref{P1} $H$ is a QHA with
antipode $S^{-1}$, co-product $\D^T$ and co-associator $\Phi^T=\Phi_{321}^{-1}$.
It follows that if $F_\d^0$ is the Drinfeld twist for this opposite structure then,
$\forall a\in H$
\begin{eqnarray*}
F_\d^0\D^T(a){(F_\d^0)}^{-1}&=&{(\D^T)}'(a)\\
&=&{(S^{-1}\ot S^{-1})}\D{(S(a))}=\D_0^T(a)
\end{eqnarray*}
since $S^{-1}$ is the antipode for this structure. On the other hand if $F_0$ is the
Drinfeld twist of equation~(\ref{E13}) we have also
\begin{equation*}
F_0^T\D^T(a){(F_0^T)}^{-1}=\D_0^T(a)
\end{equation*}
with $\D_0$ as in equation~(\ref{E8}$'$). Here we show in fact that $F_\d^0=F_0^T$.

Before proceeding we note that the Drinfeld twist is given by the canonical expression
of equation~(\ref{E9})(i) with $\gamma$ as in (\ref{E9})(ii) constructed from the
operator of (\ref{E9})(iii);~viz
\begin{eqnarray}
&& A_i\ot B_i\ot C_i\ot D_i=\left\{
\begin{array}{l}
(\Phi^{-1}\ot 1) \cdot(\D \ot 1 \ot 1) \Phi \\
\textrm{or} \\
(1 \ot \Phi)\cdot (1 \ot 1 \ot \D)\Phi^{-1}. 
\end{array} 
\right. \no 
\end{eqnarray}
This gives rise to two equivalent expansions for $\gamma$. Using the first expression
we have, in obvious notation,
\begin{eqnarray}
 A_i\ot B_i\ot C_i\ot D_i&=&{(\Phi^{-1}\ot 1)}\cdot{(\D\ot 1\ot 1)}\Phi\no\\
&=&\X_\nu X_\mu^{(1)}\ot\Y_\nu X_\mu^{(2)}\ot\Z_\nu Y_\mu\ot Z_\mu\no
\end{eqnarray}
which gives, upon substitution into (\ref{E9})(ii),
\begin{equation*}
\gamma= S{(\Y_\nu X_\mu^{(2)})}\a\Z_\nu Y_\mu\ot S{(\X_\nu X_\mu^{(1)})}\a Z_\mu
\end{equation*}
which is the expression obtained in~\cite{oqhsa}. On the other hand
using the second expression gives
\begin{eqnarray}
 A_i\ot B_i\ot C_i\ot D_i &=&{(1\ot\Phi)}\cdot{(1\ot 1\ot\D)}\Phi^{-1}\no\\
&=&\X_\mu\ot X_\nu\Y_\mu\ot Y_\nu\Z_\mu^{(1)}\ot Z_\nu\Z_\mu^{(2)}\no
\end{eqnarray}
and substituting into~(\ref{E9})(ii) gives the alternative expansion
\begin{equation}
\gamma= S{(X_\nu\Y_\mu)}\a Y_\nu\Z_\mu^{(1)}\ot S{(\X_\mu)}\a Z_\nu\Z_\mu^{(2)} \label{S4}
\end{equation}
which is equivalent to the expression above~\cite{oqhsa}.

Using (\ref{E9})(i) for the opposite structure we have for the Drinfeld twist
\begin{equation*}
F_\d^0={(S^{-1}\ot S^{-1})}\D{(X_\nu^0)}\cdot\gamma^0\cdot\D^T
{(Y_\nu^0\b^T S^{-1}{(Z_\nu^0)})}
\end{equation*}
where we have used the fact that the co-product for the opposite structure is $\D^T$,
the antipode is $S^{-1}$, with canonical elements $\a^T=S^{-1}(\a),~~~\b^T=S^{-1}(\b)$
and where we have set
\begin{equation*}
 X_\nu^0\ot Y_\nu^0\ot Z_\nu^0=\Phi^T=\Phi_{321}^{-1}~~,
\end{equation*}
which is the opposite co-associator, and where from (\ref{E9})(ii)
\begin{equation*}
\gamma^0= S^{-1}{(B_i^0)}\a^T C_i^0\ot S^{-1}{(A_i^0)}\a^T D_i^0~~
\end{equation*}
with
\begin{eqnarray}
 A_i^0\ot B_i^0\ot C_i^0\ot D_i^0 &=&{[{(\Phi^T)}^{-1}\ot 1]}\cdot
{(\D^T\ot 1\ot 1)}\Phi^T\no\\
&=&{(\Phi_{321}\ot 1)}\cdot{(\D^T\ot 1\ot 1)}\Phi_{321}^{-1}.\no
\end{eqnarray}
In obvious notation the latter is given by
\begin{equation*}
{(\Phi_{321}\ot 1)}\cdot{(\D^T\ot 1\ot 1)}\Phi_{321}^{-1}= Z_\nu\Z_\mu^{(2)}\ot
 Y_\nu\Z_\mu^{(1)}\ot X_\nu\Y_\mu\ot\X_\mu
\end{equation*}
 so that, using $\a^T=S^{-1}(\a)~~,$
 \begin{eqnarray}
 \gamma^0&=& S^{-1}{(Y_\nu\Z_\mu^{(1)})}S^{-1}(\a) X_\nu\Y_\mu\ot S^{-1}
 {(Z_\nu\Z_\mu^{(2)})}S^{-1}(\a)\X_\mu\no\\
 &\stackrel{(\ref{S4})}{=}&{(S^{-1}\ot S^{-1})}(\gamma).\no
 \end{eqnarray}
 Thus we may write, using $\b^T=S^{-1}(\b)~~,$
\begin{equation*}
 F_\d^0={(S^{-1}\ot S^{-1})}\D{(X_\nu^0)}\cdot{(S^{-1}\ot S^{-1})}\gamma
 \cdot\D^T{(Y_\nu^0 S^{-1}(\b) S^{-1}{(Z_\nu^0)})}
 \end{equation*}
 so that, substituting
 \begin{equation*}
 X_\nu^0\ot Y_\nu^0\ot Z_\nu^0=\Phi^T=\Phi_{321}^{-1}= \Z_\nu
 \ot \Y_\nu\ot \X_\nu~~,
\end{equation*}
gives
\begin{eqnarray}
F_\d^0&=& (S^{-1}\ot S^{-1})\D{(\Z_\nu)}\cdot (S^{-1}\ot S^{-1})\gamma
\cdot\D^T(\Y_\nu S^{-1}(\b)S^{-1}(\X_\nu)) \no \\
&=&(S^{-1}\ot S^{-1})\cdot [(S\ot S) \D^T (\Y_\nu S^{-1}(\X_\nu\b))\cdot
\gamma\cdot\D(\Z_\nu)] \no \\
&=&(S^{-1}\ot S^{-1}) \cdot[\D'(\X_\nu\b S(\Y_\nu)) \cdot \gamma\cdot\D
(\Z_\nu)] \no \\
&\stackrel{(\ref{E9})(i)}{=}&(S^{-1}\ot S^{-1})F_\d\stackrel{(\ref{E13})}{=}F_0^T. \no
\end{eqnarray}
Thus we have proved
\begin{Proposition}\label{P4} 
The Drinfeld twist for the opposite QHA structure of
proposition~\ref{P1} is given explicitly by
\begin{equation*}
F_\d^0={(S^{-1}\ot S^{-1})}F_\d=F_0^T.
\end{equation*}
\end{Proposition}
To see how $F_\d^T$ fits into the picture we need to consider the second Drinfeld
twist $F_0$ of Theorem~\ref{T2} associated with the co-product of equation~(\ref{E8}$'$).
We have immediately from proposition~\ref{P4}
\begin{Corollary*} The second Drinfeld twist for the opposite structure is $F_\d^T.$
\end{Corollary*}
\noindent 
\begin{proof}
Since the antipode for the opposite structure is $S^{-1}$,
Theorem~\ref{T2} implies that the second Drinfeld twist for this structure is 
$(S\ot S){(F_\d^0)}^T$ where $F_\d^0$ is the Drinfeld twist for the opposite structure, given 
explicitly in proposition~\ref{P4}. It follows that the second Drinfeld twist for the
opposite structure is
\begin{equation*}
(S \ot S) \cdot [(S^{-1} \ot S^{-1})F_\d^T]=F_\d^T.
\end{equation*}
\end{proof}
\subsection{Twisting the Drinfeld twist}
It is first useful to determine the behaviour of $\bar \g$ in equation~(\ref{E10})(ii)
under an arbitrary twist $G \in H \ot H$. Under the twisted structure induced by $G$ the
operator $\bar \g$ is twisted to $\bar \g_G$, given by equation~(\ref{E10})(ii,iii) for
the twisted structure, so that
\begin{align}
(i)&~~\bar \g_G= \bar A_i^G \b_G S(\bar D_i^G) \ot \bar B_i^G \b_G S( \bar C_i^G) \no \\
\textrm{where~~}(ii)&~~ \bar A_i^G \ot \bar B_i^G \ot \bar C_i^G \ot \bar D_i^G =
(\D_G\ot 1 \ot 1)\Phi_G^{-1}\cdot(\Phi_G \ot 1). \label{E25}
\end{align}
We have
\begin{Proposition}\label{P9}
Let $G= g_i\ot g^i \in H\ot H$ be a twist on a QHA $H$. Then
\begin{eqnarray}
\bar \g_G &=& G\cdot \D(g_i)\cdot \bar \g \cdot (S \ot S)(G^T\D^T(g^i)). \no
\end{eqnarray}
\end{Proposition}
\begin{proof}
	Throughout we write 
	\begin{eqnarray} 
	G^{-1}&=&\bar g_i \ot \bar g^i. \no
	\end{eqnarray}
	For the RHS of equation~(\ref{E25})(ii) we have
	\begin{align}
	&(\D_G\ot 1 \ot 1)\Phi_G^{-1}\cdot(\Phi_G \ot 1)=(\D_G \ot 1 \ot 1)\cdot\no \\
	&\hspace{20mm}[(1 \ot G)\cdot(1 \ot \D)G\cdot \Phi^{-1}\cdot(\D \ot 1)G^{-1}
	\cdot(G^{-1} \ot 1)]\cdot \no \\
	&\hspace{35mm}\{[(G \ot 1)\cdot(\D \ot 1)G\cdot\Phi\cdot(1 \ot \D)G^{-1}\cdot
	(1 \ot G^{-1})]\ot 1\} \no
	\end{align}
	where we have used equation~(\ref{E6}) for $\Phi_G$ and its inverse, thus
	\begin{align}
	&(\D_G\ot 1 \ot 1)\Phi_G^{-1}\cdot(\Phi_G \ot 1)=\no \\
	&\hspace{10mm}(1 \ot 1 \ot G)\cdot(\D_G \ot \D)G \cdot(\D_G \ot 1 \ot 1)\Phi^{-1}\cdot
	[(\D_G \ot 1)\D \ot 1]G^{-1}\no \\
	&\hspace{15mm}\cdot[(\D_G  \ot 1)G^{-1} \ot 1] \cdot 
	(G \ot 1 \ot 1)\cdot[(\D \ot 1)G \ot 1]\cdot(\Phi \ot 1)\no \\ 
	&\hspace{20mm}\cdot[(1 \ot \D)G^{-1}\ot 1]\cdot(1 \ot G^{-1} \ot 1) \no \\
	&\hspace{10mm}=(G \ot G)\cdot (\D \ot \D)G \cdot(\D \ot 1\ot 1)\Phi^{-1}\cdot[(\D \ot 1)\D \ot 1]G^{-1}
	 \no \\
	&\hspace{15mm}\cdot [(\D \ot 1)G^{-1}\ot 1] \cdot[(\D \ot 1)G \ot 1] \cdot (\Phi \ot 1)\no \\
	&\hspace{20mm}\cdot [(1 \ot \D)G^{-1} \ot 1]\cdot (1 \ot G^{-1} \ot 1) \no \\
	&\hspace{10mm}=(G \ot G)\cdot (\D \ot \D)G \cdot(\D \ot 1\ot 1)\Phi^{-1}\cdot[(\D \ot 1)\D \ot 1]G^{-1}
	\cdot \no \\
	&\hspace{15mm}\cdot (\Phi \ot 1)\cdot [(1 \ot \D)G^{-1} \ot 1]
	\cdot (1 \ot G^{-1} \ot 1) \no \\
	&\hspace{7mm}\stackrel{\reff{q-co}}{=}(G \ot G)\cdot (\D \ot \D)G \cdot\{(\D \ot 1\ot 1)
	\Phi^{-1}\cdot(\Phi \ot 1)\} \no \\
	&\hspace{15mm}\cdot [(1 \ot \D)\D \ot 1]G^{-1}\cdot [(1 \ot \D)G^{-1} \ot 1]\cdot (1 \ot G^{-1} \ot 1). \no
	\end{align}
	Now using the notation of equation~(\ref{E10})(iii) we have
	\begin{eqnarray}
	(\D \ot 1 \ot 1)\Phi^{-1}\cdot(\Phi \ot 1)=
	 \bar A_i \ot \bar B_i \ot \bar C_i \ot \bar D_i \no
	\end{eqnarray}
	so that in the notation of equation~(\ref{E25})(i)
	\begin{align}
	& \bar A^G_i \ot \bar B^G_i \ot \bar C^G_i \ot \bar D^G_i =
	(\D_G \ot 1 \ot 1)\Phi_G^{-1}\cdot(\Phi_G \ot 1)   \no \\
	&=(G \ot G)\cdot (\D \ot \D)G \cdot 
	\{ \bar A_i \ot \bar B_i \ot \bar C_i \ot \bar D_i\} & \no \\
	&\hspace{30mm}\cdot [(1 \ot \D)\D \ot 1]G^{-1} 
	\cdot [(1 \ot \D)G^{-1} \ot 1]\cdot (1 \ot G^{-1} \ot 1) \no \\
	&= g_sg_j^{(1)}\bar A_i\bar g_l^{(1)}\bar g_k \ot g^sg_j^{(2)}\bar B_i 
	\bar g^{(2)}_{l(1)}\bar g^k_{(1)}\bar g_m \ot g_tg^j_{(1)}\bar C_i\bar g^{(2)}_{l(2)}
	\bar g^k_{(2)}\bar g^m \ot g^tg^j_{(2)}\bar D_i \bar g^l  \no
	\end{align}
	where we have used the obvious notation, so that
	\begin{eqnarray}
	\D(g_i)&=& g_i^{(1)}\ot g_i^{(2)},\no \\
	(1\ot\D)\D(g_i)&=& g_i^{(1)}\ot \D(g_i^{(2)}) =
	 g_i^{(1)}\ot g_{i(1)}^{(2)} \ot g_{i(2)}^{(2)}, \textrm{~etc}\no
	\end{eqnarray}
	and all repeated indices are understood to be summed over.
	
	Substituting into equation~(\ref{E25})(i) gives
	\begin{eqnarray}
	\bar \g_G&=& g_sg_j^{(1)}\bar A_i\bar g_l^{(1)}\bar g_k \b_G 
	S(g^tg^j_{(2)}\bar D_i \bar g^l) \no \\
	&&\ot g^sg_j^{(2)}\bar B_i \bar g^{(2)}_{l(1)}\bar g^k_{(1)}\bar g_m \b_G 
	S(g_tg^j_{(1)}\bar C_i\bar g^{(2)}_{l(2)}\bar g^k_{(2)}\bar g^m) \no \\
	&=& g_sg_j^{(1)}\bar A_i\bar g_l^{(1)}\bar g_k \b_G 
	S(g^tg^j_{(2)}\bar D_i \bar g^l) \no \\
	&&\ot g^sg_j^{(2)}\bar B_i \bar g^{(2)}_{l(1)}\bar g^k_{(1)}\bar g_m \b_G 
	S(\bar g^m)S(\bar g^k_{(2)})S(\bar g^{(2)}_{l(2)})
	S(g_tg^j_{(1)}\bar C_i). \no
	\end{eqnarray}
	Now using
	\begin{align}
	 \bar g_m \b_G S(\bar g^m)=(\b_G)_{G^{-1}}=\b_{G^{-1}G}=\b \label{S7}
	\end{align}
	and making repeated use of equation~\reff{Sab} gives
	\begin{eqnarray}
	\bar \g_G&=& 
	 g_sg_j^{(1)}\bar A_i\bar g_l^{(1)}\bar g_k \b_G S(g^tg^j_{(2)}\bar D_i \bar g^l)\no \\
	&&\ot g^sg_j^{(2)}\bar B_i \bar g^{(2)}_{l(1)}\bar g^k_{(1)} \b 
	S(\bar g^k_{(2)})S(\bar g^{(2)}_{l(2)})
	S(g_tg^j_{(1)}\bar C_i) \no \\
	&=& g_sg_j^{(1)}\bar A_i\bar g_l \b_G S(\bar g^l)S(g^tg^j_{(2)}\bar D_i )
	\ot g^sg_j^{(2)}\bar B_i \b S(g_tg^j_{(1)}\bar C_i) \no \\
	&\stackrel{(\ref{S7})}{=}& g_sg_j^{(1)}\bar A_i \b S(\bar D_i) S(g^tg^j_{(2)})
	\ot g^sg_j^{(2)}\bar B_i \b S(\bar C_i)S(g_tg^j_{(1)}) \no \\
	&\stackrel{(\ref{E10})(ii)}{=}& (g_sg_j^{(1)} \ot g^sg_j^{(2)})\cdot \bar 
	\g \cdot (S \ot S)(g^tg^j_{(2)}
	\ot g_tg^j_{(1)}) \no \\
	&=&G \cdot \D(g_j) \cdot \bar \g \cdot (S \ot S)(G^T\cdot\D^T(g^j)) \no
	\end{eqnarray}
	\noindent which proves the result.
\end{proof}

We are now in a position to determine the action of an arbitrary twist $G \in H \ot H$
on the inverse Drinfeld twist $F_\d^{-1}$, given in equation~(\ref{E10}). Under the twisted
structure induced by $G$, $F_\d^{-1}$ is twisted to $(F^G_\d)^{-1} \equiv (F_\d^{-1})_G$,
given as in equation~(\ref{E10}), but in terms of the twisted structure, so that, with the
notation of equation~(\ref{E25}), we have from (\ref{E10})(i)
\begin{eqnarray}
(F^G_\d)^{-1}&=& \D_G (S(X_\nu^G)\a_G Y_\nu^G)\cdot \bar \g_G \cdot (S \ot S)\D^T_G(Z_\nu^G) \no
\end{eqnarray}
with $\bar \g_G$ as in proposition~\ref{P9}.

In obvious notation we may write
\begin{eqnarray}
 X_\nu^G  \ot Y_\nu^G  \ot Z_\nu^G=\Phi_G&=&(G \ot 1)\cdot(\D \ot 1)G \cdot \Phi 
\cdot (1 \ot \D)G^{-1} \cdot (1 \ot G^{-1}) \no \\
&=& g_ig_j^{(1)}X_\nu\bar g_k \ot g^ig_j^{(2)}Y_\nu\bar g^k_{(1)}\bar g_l \ot
g^jZ_\nu\bar g^k_{(2)}\bar g^l \no
\end{eqnarray}
which implies
\begin{eqnarray}
(F^G_\d)^{-1}&=&\D_G[S(g_ig_j^{(1)}X_\nu\bar g_k)\a_G g^ig_j^{(2)}Y_\nu\bar g^k_{(1)}\bar g_l]
\cdot \bar \g_G \cdot (S \ot S)\D^T_G(g^jZ_\nu\bar g^k_{(2)}\bar g^l) \no \\
&=&\D_G[S(X_\nu\bar g_k)S(g_j^{(1)})S(g_i)\a_G g^ig_j^{(2)}Y_\nu\bar
g^k_{(1)}\bar g_l]\cdot \bar \g_G \no \\
&&\cdot (S \ot S)\D^T_G(g^jZ_\nu\bar g^k_{(2)}\bar g^l). \no
\end{eqnarray}
Using
\begin{eqnarray}
 S(g_i)\a_Gg^i=(\a_G)_{G^{-1}}=\a_{G^{-1}G}=\a, \no
\end{eqnarray}
and equation~\reff{Sab}, then gives
\begin{eqnarray}
(F^G_\d)^{-1}
&=&\D_G[S(X_\nu\bar g_k)\a Y_\nu\bar g^k_{(1)}\bar g_l] 
\cdot \bar \g_G \cdot (S \ot S)\D^T_G(Z_\nu\bar g^k_{(2)}\bar g^l) \no \\
&=& G\cdot \D[S(X_\nu\bar g_k)\a Y_\nu\bar g^k_{(1)}\bar g_l]\cdot G^{-1}
\cdot \bar \g_G  \no \\
&& \cdot(S \ot S)(G^T)^{-1}\cdot (S \ot S)\D^T(Z_\nu\bar g^k_{(2)}\bar g^l)\cdot
(S \ot S)G^T \no \\
&\stackrel{prop.~(\ref{P9})}{=}& G\cdot \D[S(X_\nu\bar g_k)\a Y_\nu\bar g^k_{(1)}\bar g_l]\cdot 
\D(g_i)
\cdot \bar \g  \no \\
&&\cdot(S \ot S)\D^T (g^i)\cdot (S \ot S)\D^T(Z_\nu\bar g^k_{(2)}\bar g^l)\cdot
(S \ot S)G^T \no \\
&=& G\cdot \D[S(X_\nu\bar g_k)\a Y_\nu\bar g^k_{(1)}]\cdot 
\D(\bar g_l)\D(g_i)
\cdot \bar \g  \no \\
&&\cdot(S \ot S)\D^T (g^i)\cdot (S \ot S)\D^T (\bar g^l)\cdot
(S \ot S)\D^T(Z_\nu\bar g^k_{(2)})\no \\
&&\cdot (S \ot S)G^T \no \\
&=& G\cdot \D[S(X_\nu\bar g_k)\a Y_\nu\bar g^k_{(1)}]\cdot 
\D(\bar g_lg_i)
\cdot \g \no \\
&&\cdot (S \ot S)\D^T (\bar g^lg^i)\cdot 
(S \ot S)\D^T(Z_\nu\bar g^k_{(2)})\cdot
(S \ot S)G^T \no \\
&=& G\cdot \D[S(X_\nu\bar g_k)\a Y_\nu]\cdot 
\D(\bar g^k_{(1)})
\cdot \g  \no \\
&&\cdot (S \ot S)\D^T (\bar g^k_{(2)})\cdot 
(S \ot S)\D^T(Z_\nu)\cdot
(S \ot S)G^T \no
\end{eqnarray}
where we have used the obvious result that
\begin{eqnarray}
 \bar g_l g_i \ot \bar g^l g^i = G^{-1} G =1 \ot 1. \no
\end{eqnarray}
It then follows from proposition~\ref{P3} that
\begin{eqnarray}
(F^G_\d)^{-1}
&=& G\cdot \D[S(X_\nu)\a Y_\nu]\cdot \bar \g \cdot (S \ot S)\D^T(Z_\nu) \cdot
(S \ot S)G^T \no \\
&\stackrel{(\ref{E10})(i)}{=}&G \cdot F_\d^{-1}\cdot (S \ot S)G^T. \no
\end{eqnarray}

We have thus proved
\begin{Theorem}\label{T4}
Let $G \in H \ot H$ be a twist on a QHA $H$. Then under the twisted structure induced by
$G$, $F_\d^{-1}$ is twisted to
\begin{eqnarray}
(F^G_\d)^{-1} \equiv (F^{-1}_\d)_G=G\cdot F_\d^{-1}\cdot(S \ot S)G^T. \no
\end{eqnarray}
\end{Theorem}
Equivalently, the Drinfeld twist is twisted to 
\begin{eqnarray}
F^G_\d \equiv  (F_\d)_G=(S \ot S)(G^T)^{-1} \cdot F_\d \cdot G^{-1}. \no
\end{eqnarray}
\begin{Corollary*}
$F_0$ as in equation~(\ref{E13}) is twisted to 
\begin{eqnarray}
F^G_0 \equiv  (F_0)_G=(S^{-1} \ot S^{-1})(G^T)^{-1} \cdot F_0 \cdot G^{-1}. \no
\end{eqnarray}
\end{Corollary*}
\noindent
\begin{proof}
	Follows from the definition of $F_0 \equiv (S^{-1} \ot S^{-1})F_\d^T$ and the
	Theorem above.
\end{proof}

When $H$ is quasi-triangular the opposite structure of proposition~\ref{P1} is
obtainable, up to equivalence modulo $(S,\a,\b)$, via twisting. In such a case the
results of Section \ref{v operator} have further useful consequences. 


\section{Quasi-triangular QHAs}\label{QTQHA}

A QHA $H$ is called quasi-triangular if there exists an invertible element 
\begin{eqnarray*}
&\R&=\sum_i e_i \ot e^i \in H \ot H
\end{eqnarray*}
called the $R$-matrix, such that
\begin{eqnarray}
&(i)& \D^T(a)\R=\R\D(a),\quad \forall a\in H\no\\
&(ii)& {(\D\ot 1)}\R=\Phi_{231}^{-1} \R_{13}\Phi_{132} \R_{23}\Phi_{123}^{-1}\no\\
&(iii)& {(1\ot\D)}\R=\Phi_{312} \R_{13}\Phi_{213}^{-1} \R_{12}\Phi_{123},\label{E14}
\end{eqnarray}
where 
\begin{equation*}
\R_{12} = e_i \ot e^i \ot 1,\quad \R_{13} = e_i \ot 1 \ot e^i,\quad \text{etcetera}.
\end{equation*}

We first summarise some well known results for quasi-triangular QHAs. It was shown
in~\cite{oqhsa} that
\vskip 3mm
\noindent
{\bf Proposition \ref{P1}$~'$~~}{\it Under the opposite QHA structure of proposition~\ref{P1},
$H$ is also quasi-triangular with $R$-matrix $\R^T=T\cdot \R$, called the
opposite $R$-matrix.}
\vskip 3mm

It follows from (\ref{E14}) (ii,iii) that
\begin{equation*}
{(\e\ot 1)}\R={(1\ot \e)}\R=1
\end{equation*}
so that $\R$ qualifies as a twist. Moreover if $F\in H\ot H$ is any twist then, as
shown in~\cite{oqhsa}, $H$ is also quasi-triangular under the twisted structure of equations
(\ref{E6},~\ref{E7}) with $R$-matrix
\begin{equation}
\R_F=F^T\R F^{-1}.\label{E15}
\end{equation}
It was shown in~\cite{oqhsa} that
\begin{Proposition}\label{P5}
Under the QHA of proposition~\ref{P2}, $H$ is also quasi-triangular with $R$-matrix
\begin{equation*}
\R'=(S\ot S)\R.
\end{equation*}
\end{Proposition}

We have seen that the QHA structure of proposition~\ref{P2} is obtainable
by twisting with the Drinfeld twist $F_\d$. It was further shown in~\cite{oqhsa} that the full
structure of proposition~\ref{P5} is also obtained by twisting with $F_\d$ which,
in view of equation (\ref{E15}), is equivalent to
\begin{equation}
{(S\ot S)}\R=F_\d^T \R F_\d^{-1}.\label{E16}
\end{equation}
This result in fact follows from the following relation
\begin{equation*}
(S\ot S)\R\cdot\gamma=\gamma^T \R,
\end{equation*}
where $\gamma^T=T\cdot\gamma$, proved in~\cite{oqhsa}. In view of proposition~\ref{P3}
this last equation is equivalent to
\begin{equation*}
\R\bar\gamma=\bar\gamma^T\cdot (S\ot S) \R
\end{equation*}
where $\bar\gamma^T=T\cdot\bar\gamma$, with $\gamma$ and $\bar\gamma$ as in 
equations (\ref{E9},~\ref{E10}).

In view of (\ref{E14}) (i) the opposite co-product is obtained from $\D$ by twisting
with $\R$. In fact we have the following result proved in~\cite{oqhsa} :
\begin{Proposition}\label{P6}
The opposite structure of propositions~\ref{P1},~\ref{P1}$'$ is obtainable by twisting
with the R-matrix $\R$ but with antipode $S$ and canonical elements $\a_\R,$$~\b_\R$ respectively.
\end{Proposition}

Above $\a_\R,\b_\R$ are given by equation~(\ref{E7}), so that
\begin{align}
(i)&~~~\a_\R=m \cdot (1 \ot \a)(S \ot 1)\R^{-1}, \quad \b_\R=m \cdot (1 \ot \b)(1 \ot S)\R. \no \\
\intertext{Below we set}
(ii)&~~~\R= e_i \ot e^i,~~~\R^{-1}= \bar e_i \ot \bar e^i \no 
\intertext{in terms of which we may write}
(iii)&~~~\a_\R= S(\bar e_i)\a \bar e^i,~~~\b_\R= e_i \b S(e^i). \label{E17}
\end{align}

Thus with the co-product $\D^T$ and co-associator $\Phi^T=\Phi^{-1}_{321}$ of proposition~\ref{P1}
we have two QHA structures with differing quasi-antipodes $(S,\a_\R,\b_\R)$ and $(S^{-1},\a^T,\b^T)$
where, from proposition~\ref{P1}, $\a^T=S^{-1}(\a),~\b^T=S^{-1}(\b)$. It follows from
Theorem~\ref{T1} that
\begin{Theorem}\label{T3}
There exists a unique invertible $u \in H$ such that
\begin{equation*}
S(a)=uS^{-1}(a)u^{-1}, \textrm{~or~} S^2(a)=uau^{-1},\quad \forall a \in H
\end{equation*}
and
\begin{equation}
uS^{-1}(\a)=\a_\R,~~~\b_\R u=S^{-1}(\b). \label{E18}
\end{equation}
Explicitly,
\begin{eqnarray}
u&=& S(Y_\nu\b S(Z_\nu))\a_\R X_\nu =
 S(\Z_\nu) \a_\R \Y_\nu S^{-1}(\b)S^{-1}(\X_\nu)\no \\
u^{-1}&=&  Z_\nu \b_\R S(S(X_\nu)\a Y_\nu) =
 S^{-1}(\Z_\nu)S^{-1}(\a) \Y_\nu \b_\R S(\X_\nu). \label{E19}
\end{eqnarray}
\end{Theorem}
Above we have used the fact that the opposite QHA structure has co-associator
$\Phi^T=\Phi^{-1}_{321}$ and quasi-antipode $(S^{-1},\a^T,~\b^T)$. We have
then applied Theorem~\ref{T1} with $(\tilde S,\tilde \a,\tilde \b)=(S,\a_\R,\b_\R)$ 
to give the result.

The above gives the $u$-operator of Drinfeld-Reshetikhin~\cite{Dri90u,r}. It differs from, but is
related to, the $u$-operator of Altschuler and Coste~\cite{ac}.
To see how the latter arises, it is easily seen that $\tilde \R \equiv (\R^T)^{-1}$ also
satisfies equation~(\ref{E14}) and thus constitutes an $R$-matrix.
Thus proposition~\ref{P6} and Theorem~\ref{T2} also hold with $\R$ replaced by $\tilde \R$. 
This implies the existence of a unique invertible $\tilde u \in H$ such that
\begin{equation*}
S^2(a)=\tilde u a \tilde u^{-1},\quad \forall a \in H
\end{equation*}
and 
\begin{equation*}
\tilde u S^{-1}(\a)=\a_{\tilde \R},~~~\b_{\tilde \R} \tilde u = S^{-1}(\b)
\end{equation*}
with $\a_{\tilde \R},~\b_{\tilde \R}$ as in equation~(\ref{E17}) but with $\R$ replaced by $\tilde \R$.
Explicitly we have, in this case,
\begin{eqnarray}
\tilde u&=& S(Y_\nu\b S(Z_\nu))\a_{\tilde \R} X_\nu =
 S(\Z_\nu) \a_{\tilde \R} \Y_\nu S^{-1}(\b)S^{-1}(\X_\nu)\no \\
\tilde u^{-1}&=&  Z_\nu \b_{\tilde \R} S(S(X_\nu)\a Y_\nu) =
 S^{-1}(\Z_\nu)S^{-1}(\a) \Y_\nu \b_{\tilde \R} S(\X_\nu). \label{E20}
\end{eqnarray}
Then, as can be seen from~\cite{cas} $\tilde u$ is precisely the $u$-operator
of Altschuler and Coste.

To see the relation between $u$ and $\tilde u$ we first note that $uS(u)=S(u)u$ is central.
This follows by applying $S$ to $S(a)=uS^{-1}(a)u^{-1}$, giving
\begin{equation*}
S^2(a)=S(u^{-1})aS(u),\quad \forall a \in H.
\end{equation*}
Before proceeding it is worth noting the following 
\begin{Lemma} \label{L1}
\begin{eqnarray}
&(i)& \b_{\tilde \R}=S(u)S(\b),~~~\a_{\tilde \R}=S(\a)S(u^{-1}) \no \\
&(ii)& \b_\R=S(\tilde u)S(\b),~~~\a_\R=S(\a)S(\tilde u^{-1}). \label{E21}
\end{eqnarray}
\end{Lemma}
\begin{proof}
	By symmetry it suffices to prove (i). Now
	\begin{eqnarray}
	\b_{\tilde \R}&=& m \cdot (1 \ot \b)(1 \ot S)(\R^T)^{-1}
	= \bar e^i\b S(\bar e_i) \no \\
	&\stackrel{(\ref{E18})}{=}&\bar e^i S(\b_\R u)S(\bar e_i)
	= \bar e^iS(u)S(\b_\R)S(\bar e_i) \no \\
	&=&\bar e^iS(u)S[e_j\b S(e^j)]S(\bar e_i) \no \\
	&=&\bar e^iS(u)S^2(e^j)S(\b)S(e_j)S(\bar e_i) \no \\
	&=& S(u) S^2(\bar e^i)S^2(e^j)S(\b)S(e_j)S(\bar e_i) \no \\ 
	&=&S(u) S^2(\bar e^ie^j)S(\b)S(\bar e_ie_j)=S(u)S(\b)\no
	\end{eqnarray}
	where we have used the obvious result
	\begin{equation*}
	\bar e_ie_j \ot \bar e^ie^j=\R^{-1}\R=1 \ot 1.
	\end{equation*}
	Similarly
	\begin{eqnarray*}
	\a_{\tilde \R}&=& m \cdot (1 \ot \a)(S \ot 1)R^T
	= S(e^i)\a e_i \no \\
	&\stackrel{(\ref{E18})}{=}& S(e^i) S(u^{-1} \a_\R)e_i
	= S(e^i)S(\a_\R)S(u^{-1}) e_i \\
	&=& S(e^i)S[S (\bar e_j) \a  \bar e^j]S(u^{-1}) e_i  \\
	&=& S(e^i)S(\bar e^j)S(\a)S^2(\bar e_j)S(u^{-1}) e_i \\
	&=& S(e^i)S(\bar e^j)S(\a)S^2(\bar e_j)S^2(e_i)S(u^{-1}) \\ 
	&=& S(\bar e^je^i)S(\a)S^2(\bar e_je_i)S(u^{-1})=S(\a)S(u^{-1}).
	\end{eqnarray*} 
\end{proof}
We are now in a position to prove
\begin{Lemma}\label{L2}
\begin{equation*}
\tilde u=S(u^{-1})
\end{equation*}
\end{Lemma}
\begin{proof}
	From equation~(\ref{E20}) we have
	\begin{eqnarray}
	\tilde u &=& S(Y_\nu \b S(Z_\nu)) \a_{\tilde \R} X_\nu \no \\
	&\stackrel{(\ref{E21})(i)}{=}& S(Y_\nu \b S(Z_\nu)) S(\a)S(u^{-1}) X_\nu \no \\
	&=& S(Y_\nu \b S(Z_\nu)) S(\a)S^2(X_\nu)S(u^{-1})  \no \\
	&=&S\big[ S(X_\nu) \a Y_\nu \b S(Z_\nu)\big] S(u^{-1})  \no \\
	&\stackrel{\reff{Sphi}}{=}& S(u^{-1}). \no
	\end{eqnarray}
\end{proof}
The above result clearly shows the connection between the $u$-operator of Theorem~\ref{T2}
and that due to Altschuler and Coste. Obviously the existence of the $u$-operator in the
quasi-triangular case is a direct consequence of Theorem~\ref{T1} and proposition~\ref{P6},
the latter showing the equivalence of the opposite structure of proposition~\ref{P1} with
that due to twisting with $\R$. In the case $H$ is not quasi-triangular, this opposite
structure is not in general obtainable by a twist.

The operators $u$ and $\tilde u$ are special cases of the $v$ operator of Theorem~\ref{T1},
it follows then from Theorem~\ref{uni_of_v}, that
\begin{Theorem}\label{uni_of_u}
The operators $u$ and $\tilde u$ are invariant under twisting.
\end{Theorem}

In section~\ref{v operator} we discussed the uniqueness of the quasi-antipode $(S,\a,\b)$, 
but nothing has been said about the uniqueness of the twisted structures
or the $R$-matrix in the quasi-triangular case. This is intimately connected with the 
quasi-cocycle condition to which we now turn. 

\section{The quasi-cocycle condition}\label{QuasiCo}
The set of twists on a QHA $H$ forms a group, moreover, the
twisted structure of equations (\ref{E6}, \ref{E7}) induced on a QHA $H$ preserves this
group structure in the following sense.
\begin{Lemma}\label{L3}
Let $F,G \in H\ot H$ be twists on a QHA $H$. Then in the notation of equations~(\ref{E6},
\ref{E7})
\begin{align}
&(i)~~~\D_{FG}=(\D_G)_F,~~~\Phi_{FG}=(\Phi_G)_F \no \\
&(ii)~~~\a_{FG}=(\a_G)_F,~~~\b_{FG}=(\b_G)_F. \no \\
\intertext{Moreover, if $H$ is quasi-triangular then}
&(iii)~~~\R_{FG}=(\R_G)_F. \label{E22}
\end{align}
\end{Lemma}

In other words the structure obtained from twisting with $G$ and then with $F$ is the
same as twisting with the twist $FG$. It is important that the right hand side of
equation~(\ref{E22}) is interpreted correctly, e.g. $(\Phi_G)_F$ is given as in 
equation~(\ref{E6}) but with $\Phi$ replaced by $\Phi_G$ and $\D$ by $\D_G$ etc.

Given any QBA $H$ we may impose on a twist $F \in H \ot H$ the following condition
\begin{equation}
(F \ot 1)\cdot(\D \ot 1)F \cdot \Phi=\Phi \cdot (1\ot F)\cdot (1 \ot \D)F \label{E23}
\end{equation}
which we call the {\em quasi-cocycle condition}.

When $\Phi=1 \ot 1\ot 1$ this reduces to the usual cocycle condition on Hopf algebras.
In the notation of equation~(\ref{E6}), the quasi-cocycle condition is equivalent to
\begin{equation}
\Phi_F=\Phi. \tag{\ref{E23}'}
\end{equation}
Thus twisting on a QBA by a twist $F$ satisfying the quasi-cocycle condition results in a
QBA structure with the same co-associator.

It is thus not surprising that the quasi-cocycle condition (\ref{E23}) is intimately
related to the uniqueness of twisted structures on a QHA $H$. Indeed, if $F,G \in H\ot H$
are twists giving rise to the {\it same} QBA structure, so that
\begin{equation}
\D_F=\D_G,~~~\Phi_F=\Phi_G \label{E24}
\end{equation}
then $C\equiv F^{-1}G$ must commute with the co-product $\D$ and satisfy the quasi-cocycle
condition. Indeed in view of lemma~\ref{L3} we have
\begin{eqnarray}
\D_C=\D_{F^{-1}G}={(\D_G)}_{F^{-1}}\stackrel{(\ref{E24})}{=}&
{(\D_F)}_{F^{-1}}=\D_{F^{-1}F}=\D \no \\
\Phi_C=\Phi_{F^{-1}G}={(\Phi_G)}_{F^{-1}}\stackrel{(\ref{E24})}{=}&
{(\Phi_F)}_{F^{-1}}=\Phi_{F^{-1}F}=\Phi. \no
\end{eqnarray}
This leads to the following
\begin{Definition}
A twist $C \in H \ot H$ on any QBA $H$ is called compatible if
\begin{eqnarray}
&(i)& C \textrm{~commutes with the co-product~} \D \no \\
&(ii)& C \textrm{~satisfies the quasi-cocycle condition}~. \no
\end{eqnarray}
\end{Definition}
In other words twisting a QBA $H$ with a compatible twist $C$ gives exactly the same
QBA structure. The set of compatible twists on $H$ thus forms a subgroup of the group of twists
on $H$.

\begin{Proposition}\label{P7}
Let $F,G \in H \ot H$ be twists on a QBA $H$. Then the twisted structures induced by $F$ and
$G$ coincide if and only if there exists a compatible twist $C \in H \ot H$ such that $G=FC$.
\end{Proposition}
\begin{proof}
	We have already seen that if $F,G$ give rise to the same QBA structure then 
	$C=F^{-1}G$ is a compatible twist and $G=FC$. Conversely, suppose $C$ is a
	compatible twist and set $G=FC$. Then
	\begin{eqnarray}
	\D_G=\D_{FC}={(\D_C)}_F=\D_F \no \\
	\Phi_G=\Phi_{FC}={(\Phi_C)}_F=\Phi_F \no
	\end{eqnarray}
	so that $G$ gives precisely the same twisted structure as $F$.
\end{proof}
Setting $G=1\ot 1$ into the above gives
\begin{Corollary*}
Let $F \in H \ot H$ be a twist on a QBA $H$. Then the twisted structure induced by $F$ coincides
with the structure on $H$ if and only if $F$ is a compatible twist.
\end{Corollary*}
In view of the group properties of twists the above corollary is equivalent to
proposition~\ref{P7}.

Let $H$ be a quasi-triangular QHA with $R$-matrix $\R$ satisfying equation (\ref{E14}). 
From proposition~\ref{P6}, the opposite co-associator $\Phi^T=\Phi^{-1}_{321}$ and
co-product $\D^T$ are obtained by twisting with $\R$, so that $\Phi^T=\Phi_\R$. The proof
of this result utilises only the properties~(\ref{E14}). Hence, since
\begin{eqnarray}
\Phi=\Phi_{\R^{-1}\R}={(\Phi_\R)}_{\R^{-1}}=(\Phi^T)_{\R^{-1}} \no
\end{eqnarray}
it follows that if $Q$ is another $R$-matrix for $H$ i.e. satisfies equation~(\ref{E14}),
then we must have also
\begin{eqnarray}
{(\Phi^T)}_{Q^{-1}}=\Phi.\no
\end{eqnarray}
Then $Q^{-1}\R$ must qualify as a compatible twist. Indeed it obviously 
commutes with $\D$, while as to the quasi-cocycle condition, we have 
\begin{eqnarray} 
\Phi_{Q^{-1}\R}={(\Phi_\R)}_{Q^{-1}}=(\Phi^T)_{Q^{-1}}=\Phi. \no
\end{eqnarray}
Note that $(Q^T)^{-1}, (\R^T)^{-1}$ also determine $R$-matrices so the following
must all determine compatible twists: 
$Q^{-1}\R,~Q^T\R,~\R^{-1}Q,~\R^T Q$.
In particular $\R^T\R$ must determine a compatible twist, as may be verified directly. 

With the notation of section \ref{twisting}, it is easily seen that the operator
\begin{equation}
A=\D(u^{-1})F_\delta^{-1}(u \ot u)F_0=F_\delta^{-1}(u \ot u)F_0\D(u^{-1}) \label{E24AC}
\end{equation}
commutes with $\D$. This operator appears in the work of Altschuler and Coste~\cite{ac} in
connection with ribbon QHAs. The operator $A$ satisfies the quasi-cocycle condition and thus 
determines a compatible twist.

For general QBAs $H$, to see that there are sufficiently many compatible twists, we have
\begin{Lemma}\label{L4}
Let $z \in H$ be an invertible central element. Then
\begin{eqnarray}
C=(z\ot z)\D(z^{-1}) \no
\end{eqnarray}
is a compatible twist.
\end{Lemma}
\begin{proof}
	Obviously $C$ commutes with the co-product $\D$ so it remains to prove that it satisfies
	the quasi-cocycle condition. To this end note that
	\begin{eqnarray}
	(C \ot 1)(\D \ot 1)C &=&
	(z \ot z \ot 1)(\D(z^{-1}) \ot 1)(\D(z)\ot z)(\D \ot 1)\D(z^{-1}) \no \\
	&=&(z \ot z \ot z)(\D \ot 1)\D(z^{-1}) \label{S5}
	\end{eqnarray}
	and similarly
	\begin{eqnarray}
	(1 \ot C)(1 \ot \D)C&=&
	(1 \ot z \ot z)(1 \ot \D(z^{-1}))(z \ot \D(z)) (1 \ot \D)\D(z^{-1}) \no \\
	&=&(z \ot z \ot z)(1 \ot \D)\D(z^{-1}) \label{S6}
	\end{eqnarray}
	thus
	\begin{eqnarray}
	(C \ot 1)(\D \ot 1)C\Phi&\stackrel{~~(\ref{S5})~~}{=}&
	(z \ot z \ot z)(\D \ot 1)\D(z^{-1}) \Phi \no \\
	&\stackrel{\reff{q-co}~~}{=}&(z \ot z \ot z) \Phi (1\ot\D)\D(z^{-1})\no \\
	&\stackrel{~~(\ref{S6})~~}{=}&(z \ot z \ot z) \Phi
	 (z^{-1}\ot z^{-1}\ot z^{-1}) \no \\
	&& (1 \ot C)(1 \ot \D)C \no \\
	&=&\Phi(1 \ot C)(1 \ot \D)C. \no
	\end{eqnarray} 
\end{proof}
With $C$ as in the lemma, we see that
\begin{eqnarray}
(\e \ot 1)C=(1 \ot \e)C=\e(z). \no
\end{eqnarray}
Thus, strictly speaking, $\e(z^{-1})C$ qualifies as a compatible twist. 

Following Altschuler and Coste~\cite{ac}, a quasi-triangular QHA is called a ribbon QHA
if the operator $A$ of equation~(\ref{E24AC}) is given by
\begin{eqnarray}
A=(v \ot v)\D(v^{-1}) \no
\end{eqnarray}
for a certain invertible central element $v$, related to the $u$-operator $u$. This is 
consistent with the lemma above and the fact that $A$
determines a compatible twist.

In the case of ribbon Hopf algebras, we have $\R^T\R=(v \ot v)\D(v^{-1})$, so
that the compatible twist $\R^T\R$ is also of the form of lemma~\ref{L4}. This
may not be the case for quasi-triangular QHAs in general.

It is worth noting that if $H$ is a QHA and $C \in H \ot H$ a compatible twist then
$H$ is also a QHA under the twisted structure induced by $C$ with exactly the same 
co-product $\D$, co-unit $\e$, co-associator $\Phi$, antipode $S$, but with canonical
elements given by equation~(\ref{E7}); viz
\begin{eqnarray}
\a_C=m\cdot(S \ot 1)(1 \ot \a)C^{-1},~~~\b_C=m\cdot(1 \ot S)(1\ot \b)C. \no
\end{eqnarray}
In view of Theorem~\ref{T1} and it corollary, we have immediately
\begin{Proposition}\label{P8}
Suppose $C \in H \ot H$ is a compatible twist on a QHA $H$. Then there exists a unique
invertible central element $z \in H$ such that
\begin{eqnarray}
z \a=\a_C,~~~\b_Cz=\b. \no
\end{eqnarray}
Explicitly
\begin{eqnarray}
z&=& S(X_\nu)\a_C Y_\nu \b S(Z_\nu)= \X_\nu \b S(\Y_\nu)\a_C \Z_\nu \no \\
z^{-1}&=& S(X_\nu)\a Y_\nu \b_C S(Z_\nu)= \X_\nu \b_C S(\Y_\nu) \a \Z_\nu. \no
\end{eqnarray}
\end{Proposition}

In the case $H$ is quasi-triangular we have seen that $C=\R^T\R$ is a compatible twist.
Since the latter form a group we have the infinite family of compatible twists 
$C=(\R^T\R)^m,~m \in \mathbb Z$, in which case the central elements $z^{\pm 1}$ of
proposition~\ref{P8} give the quadratic invariants of~\cite{cas}.

We conclude this section by noting, in the quasi-triangular case, that twisting the
Drinfeld twist with the $R$-matrix $\R$ gives, from Theorem~\ref{T4}, the twisted
Drinfeld twist
\begin{eqnarray}
F_\d^\R \equiv (F_\d)_\R=(S \ot S)(\R^T)^{-1}\cdot F_\d \cdot \R^{-1}. \no
\end{eqnarray}
On the other hand, since $(\R^T)^{-1}$ is an $R$-matrix we have, from eq.~(\ref{E16}),
\begin{eqnarray}
(S \ot S)(\R^T)^{-1}=F_\d^T(\R^T)^{-1}F_\d^{-1} \no
\end{eqnarray}
which implies
\begin{eqnarray}
F^\R_\d=F^T_\d(\R^T)^{-1} \cdot \R^{-1}=F^T_\d(\R\R^T)^{-1} \no
\end{eqnarray}
where $\R\R^T$ and its inverse are compatible twists under the opposite structure. This shows
that $F_\d^T$ will give rise to a Drinfeld twist under the opposite structure of 
proposition~\ref{P6} induced by twisting with $\R$ (which has antipode $S$ rather than
$S^{-1}$). Applying $T$ to the equation above gives
\begin{eqnarray}
(F^R_\d)^T=F_\d(\R^T\R)^{-1} \no
\end{eqnarray}
which shows that, since $\R^T\R$ and its inverse are compatible twists, $(F^\R_\d)^T$ also
gives rise to a Drinfeld twist on $H$.


\section{Quasi-dynamical QYBE}
Throughout we assume $H$ is a quasi-triangular QHA with $R$-matrix $\R$ satisfying (\ref{E14})
which we reproduce here:
\begin{align}
(i)&~~\D^T(a)\R=\R \D(a),\quad \forall a \in H \no \\
(ii)&~~(\D \ot 1)\R=\Phi^{-1}_{231} \R_{13} \Phi_{132} \R_{23} \Phi^{-1}_{123} \no \\
(iii)&~~(1 \ot \D)\R=\Phi_{312} \R_{13} \Phi^{-1}_{213} \R_{12} \Phi_{123}. \tag{\ref{E14}$'$}\label{E14'}
\end{align}
Applying $T \ot 1$ to (ii) and $1 \ot T$ to (iii) then gives
\begin{align}
(ii')&~~(\D^T \ot 1)\R=\Phi^{-1}_{321} \R_{23} \Phi_{312} \R_{13} \Phi^{-1}_{213} \no \\
(iii')&~~(1 \ot \D^T)\R=\Phi_{321} \R_{12} \Phi^{-1}_{231} \R_{13} \Phi_{132}. \no
\end{align}
It follows that 
\begin{equation*}
\R_{12}(\D \ot 1)\R=(\D^T \ot 1)\R \cdot \R_{12}
\end{equation*}
from which we deduce that $\R$ must satisfy the quasi-QYBE:
\begin{equation}
\R_{12}\Phi^{-1}_{231} \R_{13} \Phi_{132} \R_{23} \Phi^{-1}_{123}=
\Phi^{-1}_{321} \R_{23} \Phi_{312} \R_{13} \Phi^{-1}_{213}\R_{12}. \label{E42}
\end{equation}

If we twist $H$ with a twist $F \in H\ot H$ then $H$ is also a quasi-triangular QHA under the twisted
structure (\ref{E6},\ref{E7}) induced by $F$ with universal $R$-matrix
\begin{equation*}
\R_F=F^T \R F^{-1}.
\end{equation*}

Following equation (\ref{E23}) we say a twist $F(\lambda) \in H\ot H$ satisfies the
shifted quasi-cocycle condition if 
\begin{equation}
[F(\lambda) \ot 1]\cdot(\D \ot 1)F(\lambda)\cdot \Phi= \Phi\cdot[1 \ot F(\lambda+h^{(1)})]
\cdot (1 \ot \D)F(\lambda) \label{E43}
\end{equation}
where $\lambda \in H$ depends on one (or possibly several) parameters and $h \in H$ is fixed.
Alternatively, we may write in obvious notation
\begin{align}
F_{12}(\lambda)\cdot(\D \ot 1)F(\lambda)\cdot \Phi= \Phi\cdot F_{23}(\lambda+h^{(1)})
\cdot (1 \ot \D)F(\lambda). \tag{\ref{E43}$'$}\label{E43'}
\end{align}
When $h=0$, this reduces to the quasi-cocycle condition (\ref{E23}) satisfied by $F=F(\lambda)$.
When $\Phi =1 \ot 1\ot 1$ (i.e. the normal Hopf-algebra case) equation (\ref{E43}) reduces
to the usual shifted cocycle condition.

Twisting $H$ with a twist $F$ satisfying the (unshifted) quasi-cocycle condition results in a
QHA with the same co-associator $\Phi$, co-unit $\e$ and antipode $S$ but with the twisted
co-product $\D_F$, $R$-matrix $\R_F$ (and canonical elements $\a_F,\b_F$).
We now consider twisting $H$ with a twist $F=F(\lambda)$ satisfying the shifted
condition (\ref{E43}). Then under this twisted structure $H$ is also a quasi-triangular QHA with
the same co-unit $\e$ and antipode $S$ but with the co-associator $\Phi(\lambda)=\Phi_{F(\lambda)}$,
and co-product and $R$-matrix given by
\begin{equation}
\D_\lambda(a)=F(\lambda)\D(a)F(\lambda)^{-1},\quad \forall a \in H,
~~~\R(\lambda)=F^T(\lambda)\R F(\lambda)^{-1} \label{E44}
\end{equation}
with canonical elements $\a_\lambda=\a_{F(\lambda)},\b_\lambda=\b_{F(\lambda)}$.

In view of equation (\ref{E43}$'$) we have for the co-associator
\begin{eqnarray}
\Phi(\lambda)&=&F_{12}(\lambda)\cdot(\D \ot 1)F(\lambda)\cdot \Phi 
\cdot (1 \ot \D)F(\lambda)^{-1}\cdot F_{23}(\lambda)^{-1} \no \\
&=&\Phi \cdot F_{23}(\lambda+h^{(1)}) \cdot (1 \ot \D)F(\lambda) \cdot(1 \ot \D)F(\lambda)^{-1}  
\cdot F_{23}(\lambda)^{-1} \no \\
&=&\Phi \cdot F_{23}(\lambda+h^{(1)}) \cdot F_{23}(\lambda)^{-1} \label{E45}
\end{eqnarray}
which implies
$$\Phi(\lambda)^{-1}= F_{23}(\lambda) \cdot F_{23}(\lambda+h^{(1)})^{-1} \cdot \Phi^{-1}.$$
In the Hopf-algebra case equation (\ref{E45}) reduces to the expression for
$\Phi(\lambda)$ obtained in~\cite{oqhsa} ($\Phi= 1\ot 1\ot1)$.

Under the above twisted structure equation (\ref{E14})(ii) becomes
$$(\D_\lambda \ot 1)\R(\lambda)=\Phi_{231}(\lambda)^{-1} \cdot\R_{13}(\lambda) \cdot \Phi_{132}(\lambda)
\cdot \R_{23}(\lambda) \cdot \Phi^{-1}_{123}(\lambda).$$ 
Now
\begin{align}
\Phi_{132}(\lambda)&=(1\ot T)\Phi_{123}(\lambda) \no \\
&\stackrel{(\ref{E45})}{=}\Phi_{132} \cdot
F^T_{23}(\lambda+h^{(1)}) \cdot F^T_{23}(\lambda)^{-1} \label{S14}
\end{align}
which implies
\begin{eqnarray*}
(\D_\lambda \ot 1)\R(\lambda)&=&\Phi_{231}(\lambda)^{-1} \cdot \R_{13}(\lambda) \cdot \Phi_{132}
\cdot F^T_{23}(\lambda+h^{(1)}) \\
&& \cdot F^T_{23}(\lambda)^{-1} \cdot \R_{23}(\lambda) \cdot \Phi^{-1}_{123}(\lambda) \no \\
&\stackrel{(\ref{E45})}{=}&\Phi_{231}(\lambda)^{-1} \cdot \R_{13}(\lambda) \cdot \Phi_{132}
\cdot F^T_{23}(\lambda+h^{(1)})  \\
&&\cdot F^T_{23}(\lambda)^{-1}\cdot\R_{23}(\lambda) \cdot F_{23}(\lambda) 
\cdot F_{23}(\lambda+h^{(1)})^{-1} \cdot \Phi^{-1}_{123} \no \\
&\stackrel{(\ref{E44})}{=}&\Phi_{231}(\lambda)^{-1} \cdot \R_{13}(\lambda) \cdot \Phi_{132}
\cdot \R_{23}(\lambda+h^{(1)}) \cdot \Phi^{-1}_{123}. \no
\end{eqnarray*}

Similarly equation (\ref{E14})(iii) becomes
$$(1 \ot \D_\lambda)\R(\lambda)=\Phi_{312}(\lambda) \cdot \R_{13}(\lambda) \cdot \Phi^{-1}_{213}(\lambda)
\cdot \R_{12}(\lambda) \cdot \Phi_{123}(\lambda).$$
Now
\begin{align}
\Phi_{312}(\lambda)&=(T \ot 1)(1 \ot T)\Phi_{123}(\lambda) \no \\
&\stackrel{(\ref{E45})}{=}(T \ot 1)[\Phi_{132}\cdot F^T_{23}(\lambda+h^{(1)}) \cdot F^T_{23}(\lambda)^{-1}] \no \\
&=\Phi_{312}\cdot F^T_{13}(\lambda+h^{(2)}) \cdot F^T_{13}(\lambda)^{-1} \no
\end{align}
while
\begin{align}
\Phi_{213}^{-1}(\lambda)&=(T \ot 1)\Phi(\lambda)^{-1} \no \\
&\stackrel{(\ref{E45})}{=}(T \ot 1)[F_{23}(\lambda)\cdot F_{23}(\lambda+h^{(1)})^{-1} \cdot \Phi^{-1}] \no \\
&=F_{13}(\lambda) \cdot F_{13}(\lambda+h^{(2)})^{-1} \cdot \Phi_{213}^{-1}. \no
\end{align}
Therefore
\begin{align}
(1 \ot \D_\lambda)\R(\lambda)&= \Phi_{312} \cdot F^T_{13}(\lambda+h^{(2)}) \cdot F^T_{13}(\lambda)^{-1}
\cdot \R_{13}(\lambda)  \no \\
&\hspace{20mm}\cdot F_{13}(\lambda) \cdot F_{13}(\lambda+h^{(2)})^{-1}\cdot\Phi^{-1}_{213}
\cdot \R_{12}(\lambda) \cdot \Phi_{123}(\lambda) \no \\
&\stackrel{(\ref{E44})}{=}  \Phi_{312} \cdot R_{13}(\lambda+h^{(2)})\cdot\Phi^{-1}_{213}
\cdot R_{12}(\lambda) \cdot \Phi_{123}(\lambda). \no
\end{align}
We thus arrive at
\begin{Lemma}\label{L14} $\R(\lambda)$ satisfies the co-product properties
\begin{eqnarray}
&(i)& (\D_\lambda \ot 1)\R(\lambda)=\Phi_{231}^{-1}(\lambda) \cdot \R_{13}(\lambda) \cdot \Phi_{132}
\cdot \R_{23}(\lambda+h^{(1)}) \cdot \Phi^{-1}_{123} \no \\
&(ii)& (1 \ot \D_\lambda)\R(\lambda)=\Phi_{312} \cdot R_{13}(\lambda+h^{(2)})\cdot\Phi^{-1}_{213}
\cdot R_{12}(\lambda) \cdot \Phi_{123}(\lambda) \no \\
&(iii)& (\D^T_\lambda \ot 1)\R(\lambda)=\Phi_{321}^{-1}(\lambda) \cdot \R_{23}(\lambda) \cdot \Phi_{312}
\cdot \R_{13}(\lambda+h^{(2)}) \cdot \Phi^{-1}_{213} \no \\
&(iv)& (1 \ot \D^T_\lambda)\R(\lambda)=\Phi_{321} \cdot R_{12}(\lambda+h^{(3)})\cdot\Phi^{-1}_{231}
\cdot R_{13}(\lambda) \cdot \Phi_{132}(\lambda). \no \\
&&\label{E46}
\end{eqnarray}
\end{Lemma}
\begin{proof}
	We have already proved (i) and (ii) while (iii) follows by applying $(T \ot 1)$ to (i) and
	(iv) by applying $(1 \ot T)$ to (ii).
\end{proof}

We are now in a position to determine the QQYBE (\ref{E42}) satisfied by $\R=\R(\lambda)$ for this twisted
structure. We have
\begin{eqnarray}
&\mbox{}&\R_{23}(\lambda) \cdot \Phi_{312} \cdot \R_{13}(\lambda+h^{(2)}) \cdot \Phi^{-1}_{213} 
\cdot R_{12}(\lambda) \no \\
&\stackrel{(\ref{E46})(ii)}{=}&\R_{23}(\lambda) \cdot (1 \ot \D_\lambda)\R(\lambda) 
\cdot \Phi^{-1}_{123}(\lambda) \no \\
&\stackrel{(\ref{E14})(i)}{=}&(1 \ot \D^T_\lambda)\R(\lambda) \cdot \R_{23}(\lambda)
\cdot \Phi^{-1}_{123}(\lambda) \no \\
&\stackrel{(\ref{E46})(iv)}{=}&\Phi_{321} \cdot R_{12}(\lambda+h^{(3)})\cdot\Phi^{-1}_{231}
\cdot R_{13}(\lambda) \cdot \Phi_{132}(\lambda) \cdot \R_{23}(\lambda)
\cdot \Phi^{-1}_{123}(\lambda) \no
\end{eqnarray}
where for the last three terms we have
\begin{eqnarray*}
\Phi_{132}(\lambda) \cdot \R_{23}(\lambda) \cdot \Phi^{-1}_{123}(\lambda) 
&\stackrel{(\ref{E45},\ref{S14})}{=}& \Phi_{132} \cdot F^T_{23}(\lambda+h^{(1)})
\cdot F^T_{23}(\lambda)^{-1} \\
&& \cdot R_{23}(\lambda) \cdot F_{23}(\lambda) \cdot F_{23}(\lambda+h^{(1)})^{-1} \cdot \Phi^{-1}_{123} \no \\
&\stackrel{(\ref{E44})}{=}& \Phi_{132} \cdot \R_{23}(\lambda+h^{(1)}) \cdot \Phi^{-1}_{123}. \no
\end{eqnarray*}
Hence
\begin{eqnarray}
&\mbox{}&\R_{23}(\lambda) \cdot \Phi_{312} \cdot \R_{13}(\lambda+h^{(2)}) \cdot \Phi^{-1}_{213} 
\cdot R_{12}(\lambda) \no \\
&=&\Phi_{321}\cdot \R_{12}(\lambda+h^{(3)}) \cdot \Phi^{-1}_{231} \cdot 
\R_{13}(\lambda) \cdot \Phi_{132} \cdot \R_{23}(\lambda+h^{(1)}) \cdot \Phi^{-1}_{123}. \no
\end{eqnarray}
We thus arrive at 
\begin{Proposition}\label{P19}
$\R(\lambda)$ satisfies the quasi-dynamical QYBE
\begin{eqnarray}
&\mbox{}&\R_{12}(\lambda+h^{(3)}) \cdot \Phi^{-1}_{231} \cdot \R_{13}(\lambda) \cdot \Phi_{132}
\cdot \R_{23}(\lambda+h^{(1)}) \cdot \Phi^{-1}_{123} \no \\
&=& \Phi^{-1}_{321} \cdot \R_{23}(\lambda) \cdot \Phi_{312} \cdot \R_{13}(\lambda+h^{(2)}) \cdot \Phi^{-1}_{213} 
\cdot \R_{12}(\lambda). \label{E47} 
\end{eqnarray}
\end{Proposition}
In the Hopf algebra case ($\Phi =1 \ot 1 \ot 1$) equation (\ref{E47}) reduces to the usual
dynamical QYBE. If we set $h=0$ then equation (\ref{E47}) reduces to the quasi-QYBE
(\ref{E42}) satisfied by $\R=\R(\lambda)$. Hence the term quasi-dynamical QYBE for (\ref{E47}):
we could, alternatively, refer to (\ref{E47}) as the dynamical quasi-QYBE (dynamical QQYBE),
since it is obviously the quasi-Hopf algebra analogue of the usual dynamical QYBE. 


With respect to the QHA structure of propositions \ref{P2}, \ref{P2}$'$ we have the $R$-matrices
$$\R'(\lambda)=(S \ot S)\R(\lambda),~~~\R_0(\lambda)=(S^{-1} \ot S^{-1})\R(\lambda)$$
respectively. Then applying $(S \ot S \ot S),~(S^{-1} \ot S^{-1} \ot S^{-1})$ respectively to
equation (\ref{E47}) it follows that both of these $R$-matrices satisfy the opposite
quasi-dynamical QYBE
\begin{eqnarray}
&\mbox{}& \tilde \R_{12}(\lambda)\cdot \tilde\Phi^{-1}_{231} \cdot \tilde\R_{13}(\lambda+h^{(2)}) \cdot 
\tilde\Phi_{132} \cdot \tilde\R_{23}(\lambda) \cdot \tilde\Phi^{-1}_{123} \no \\
&=& \tilde\Phi^{-1}_{321} \cdot \tilde\R_{23}(\lambda+h^{(1)}) \cdot \tilde\Phi_{312} \cdot \tilde\R_{13}(\lambda)
\cdot \tilde\Phi^{-1}_{213} \cdot \tilde\R_{12}(\lambda+h^{(3)}) \no
\end{eqnarray}
where $\tilde\Phi$ is the co-associator of propositions \ref{P2},\ref{P2}$'$ and 
$\tilde\R(\lambda)$ denotes $\R'(\lambda)$, $\R_0(\lambda)$ respectively. Moreover applying $(T \ot 1)((1 \ot T)(T \ot 1)$ to equation
(\ref{E47}) it is easily seen that $\R^T(\lambda)$ also satisfies the above opposite quasi-dynamical QYBE
but with respect to the opposite co-associator $\Phi^T$ of proposition \ref{P1}. 

We anticipate that the quasi-dynamical QYBE will play an important role in obtaining elliptic solutions
to the QQYBE from trigonometric ones via twisted QUEs. Of particular interest is the quasi-dynamical QYBE
for elliptic quantum groups.

\newpage

\end{document}